\documentclass[11pt]{article}
\usepackage{amsmath,amsthm,amssymb}
\usepackage{thmtools,thm-restate}
\usepackage[noadjust]{cite}
\usepackage{graphicx}
\usepackage{mathrsfs}
\usepackage{hyperref}
\usepackage{mathtools}
\hypersetup{colorlinks=true,linkcolor=blue,filecolor=blue,citecolor=blue}
\usepackage{cleveref}
\usepackage{cite}
\usepackage{makecell}
\usepackage{array}
\usepackage[defaultlines=3,all]{nowidow}
\usepackage{microtype}
\usepackage{multirow}
\usepackage{mdwtab}
\usepackage{comment}
\newtheorem{theorem}{Theorem}[section]
\newtheorem*{theorem*}{Theorem}
\newtheorem{lemma}[theorem]{Lemma}
\newtheorem{corollary}[theorem]{Corollary}
\newtheorem{conjecture}[theorem]{Conjecture}
\newtheorem*{conjecture*}{Conjecture}

\newtheorem{proposition}[theorem]{Proposition}
\newtheorem{problem}{Problem}

\newtheorem{claim}{Claim}
\newtheorem{observation}{Observation}

\DeclareMathOperator{\g}{g}
\DeclareMathOperator{\mad}{mad}
\newenvironment{clproof}{ \trivlist
        \item[\hskip\labelsep
        \emph{Proof of the claim}.]\ignorespaces
}{\hfill$\vartriangleleft$\medskip
        
}

\begin{document}
\title{\bf  Semistrong edge colorings of planar graphs
}
\author{Yuquan Lin
\  and Wensong Lin\footnote{Corresponding author. E-mail address: wslin@seu.edu.cn}\\
{\small School of Mathematics, Southeast University, Nanjing 210096, P.R. China}}
\date{}
\maketitle
\vspace*{-1cm}
\begin{abstract}
Strengthened notions of a matching $M$ of a graph $G$ have been considered, requiring that the matching $M$ has some properties with respect to the subgraph $G_M$ of $G$ induced by the vertices covered by $M$:
If $M$ is the unique perfect matching of $G_M$, then $M$ is a \emph{uniquely restricted matching} of $G$; if all the edges of $M$ are pendant edges of $G_M$, then $M$ is a \emph{semistrong matching} of $G$; if all the vertices of $G_M$ are pendant, then $M$ is an \emph{induced matching} of $G$.
Strengthened notions of edge coloring and of the chromatic index follow.

In this paper, we consider the maximum semistrong chromatic index of planar graphs with given maximum degree $\Delta$.
We prove that graphs with maximum average degree less than ${14}/{5}$ have semistrong chromatic index (hence uniquely restricted chromatic index) at most $2\Delta+4$, and we reduce the bound to $2\Delta+2$ if the maximum average degree is less than ${8}/{3}$.
These cases cover, in particular,  the cases of planar graphs with girth at least 7 (resp. at least 8). 

Our result makes some progress on the conjecture of 
Lu{\v{z}}ar, Mockov{\v{c}}iakov{\'a} and Sot{\'a}k [J.~Graph Theory 105 (2024) 612--632], which asserts that every planar graph $G$  has a semistrong edge coloring with $2\Delta+C$ colors, for some universal constant $C$.
(Note that such a conjecture would fail for strong edge coloring as there exist graphs with arbitrarily large maximum degree that are not strongly $(4\Delta-5)$-edge-colorable.) We provide an example of a planar graph showing that the maximum semistrong chromatic index of planar graphs with maximum degree $\Delta$ is at least $2\Delta+4$.
\end{abstract}




 \section{Introduction}


We consider  finite undirected simple graphs in this paper.
   Let $G=(V(G),E(G))$ be a  graph, and let  $\Delta$ be its  maximum degree. A vertex of a graph $G$ is \emph{pendant} if it has degree $1$; an edge of $G$ is \emph{pendant} if it is incident to a pendant vertex.
A matching $M$ of a graph $G$ is a set of pairwise disjoint edges, and we denote by $G_M$  the subgraph of $G$ induced by the vertices covered by $M$. (Note that $M$ is a perfect matching of $G_M$.)
A \emph{$\mathscr{P}$-matching} \cite{GHHL2005} of a graph $G$ is a matching $M$ such that $G_M$ has the property $\mathscr{P}$.
From the weakest to the strongest property $\mathscr P$, 
we consider the following {$\mathscr{P}$-matching} (See \Cref{fig:matching}):
\begin{itemize}
    \item \emph{uniquely restricted matching}: $G_M$ has a unique perfect matching,
    \item \emph{semistrong matching}: every edge of $M$ is pendant in $G_M$,
    \item \emph{induced matching} (or \emph{strong matching}): every vertex of $G_M$ is pendant.
\end{itemize}

 
  \begin{figure}[ht]
    \centering
    {\includegraphics[width=\textwidth]{matching.pdf}}
    \caption{Different types of matching (heavy edges), from weaker to stronger.}
    \label{fig:matching}
\end{figure}

\begin{figure}[h!t]
	\begin{center}
		\includegraphics[width=.8\textwidth]{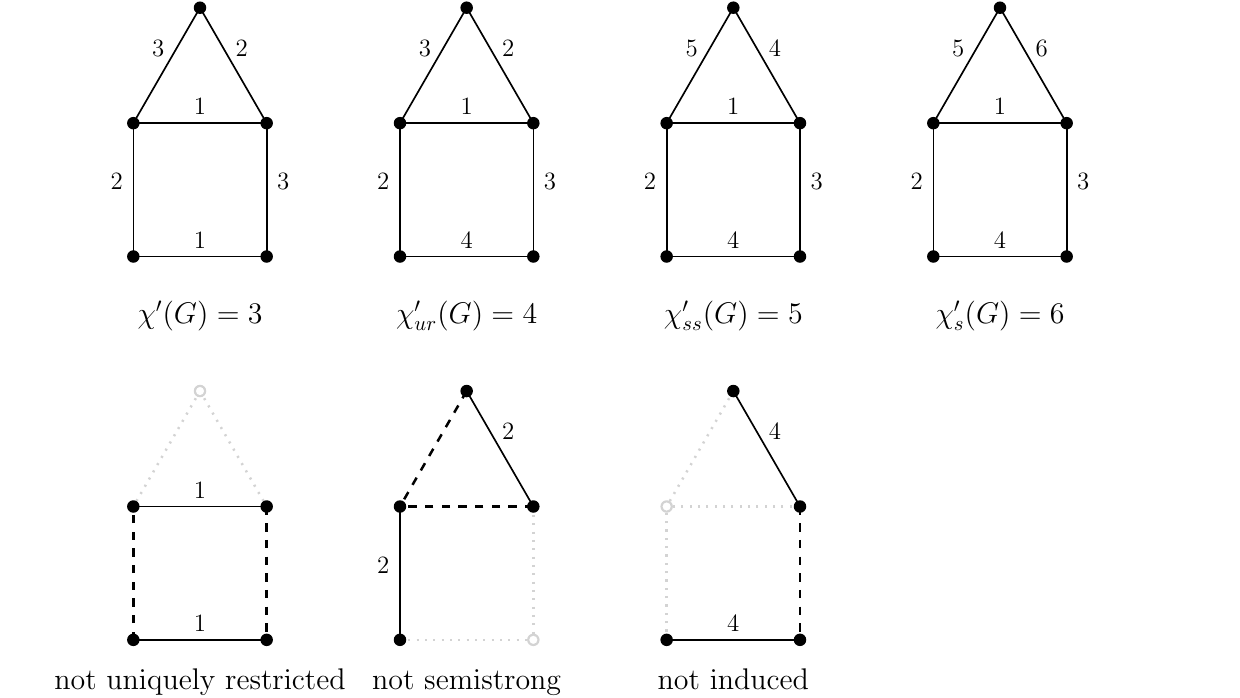}
		\caption{Different types of edge coloring (proper, uniquely restricted, semistrong, strong).}
		\label{fig:house}
	\end{center}
\end{figure}

 Each type of matching derives a corresponding edge coloring.
 An edge coloring of a graph $G$ is an assignment of colors to its edges such that each color class is a matching.
The  \emph{uniquely restricted edge coloring} \cite{BRS2017},  \emph{semistrong edge coloring} \cite{GH2005}, and  \emph{strong edge coloring} \cite{FJ1983} of a graph $G$ is an edge coloring such that each color class is a uniquely restricted matching, a semistrong matching, and an induced matching, respectively (See \Cref{fig:house}).

Notice that a matching $M$ is uniquely restricted if and only if there is no $M$‐alternating cycle (a cycle that alternates between edges in $M$ and edges not in $M$ ) in $G$ \cite{GHL2001}. 
The uniquely restricted edge coloring is a strengthening of the   \emph{acyclic edge coloring} introduced in \cite{F1978}, which is an edge coloring of a graph $G$ such that each cycle of $G$ gets at least 3 colors.
The \emph{chromatic index} $\chi'(G)$, \emph{acyclic  chromatic index} $a'(G)$,  \emph{uniquely restricted  chromatic index} $\chi_{ur}'(G)$,  \emph{semistrong  chromatic index} $\chi_{ss}'(G)$, and  \emph{strong  chromatic index}  $\chi_{s}'(G)$ of a graph $G$ are defined as the minimum number of colors required for the corresponding colorings.
Given the strength of the constraints corresponding to these colorings, every graph $G$ satisfies
$$\chi'(G)\le a'(G)\le \chi_{ur}'(G)\le \chi_{ss}'(G)\le \chi_{s}'(G).$$

Both the strong edge coloring and the acyclic edge coloring have been widely studied \cite{B2010,DYZ2019}.  
Remark that all  known results for strong edge coloring  are applicable to  semistrong edge coloring. 
In the planar case,   the following upper bounds have been obtained (See Table \ref{tab:bound}). (In the table, $\g(G)$ denotes  the \emph{girth} of the graph  $G$,  that is the length of a shortest cycle of $G$.)

On the one hand,
 Faudree, Schelp, Gy{\'a}rf{\'a}s,  and Tuza \cite{FSGT1990} proved that   $\chi'_{s}(G)\le  4\Delta + 4$ for any planar graph $G$, and
the coefficient $4$  is  best possible.
Indeed, they  constructed a class of planar graphs with maximum degree $\Delta\ge 2$ and strong chromatic index $4\Delta-4$.
Also, the  \emph{triangular prism} is a 3-regular planar graph with  semistrong chromatic index 9 ($=4\Delta-3$) \cite{FSGT1990} and there  exist  two planar graphs of maximum degree $4$ with semistrong chromatic index 13 ($=4\Delta-3$) \cite{WSWC2018}.
On the other hand, Shu and Lin showed very recently \cite{SL2024} that $a'(G)\le\Delta+5$ and this is currently the best known upper bound for planar graphs.
However,  the Acyclic Edge Coloring  Conjecture (AECC) states that $a'(G)\le \Delta+2$ holds	for any graph $G$.  This conjecture was independently proposed by Fiamčik \cite{F1978} and by Alon, Sudakov, and  Zaks \cite{ASZ2001}.

In contrast, the uniquely restricted edge coloring and the semistrong edge coloring have been relatively little studied on planar graphs.
This motivates us to study the upper bounds of $\chi_{ur}'(G)$ and  $\chi_{ss}'(G)$ on planar graphs.
For a study of these edge colorings (for general graphs), we refer the reader to \cite{BRS2019,BRS2017} and \cite{GH2005,LMS2024,LL2023}, respectively.

\medskip

The work presented in this paper is motivated by the following conjecture.

\begin{restatable}[{\cite[Conjecture 3]{LMS2024}}]{conjecture}{LMSConj}
\label{conj-planar-semi-L} 
 	There is a (small) constant $C$ such that for any planar graph $G$ with  maximum degree $\Delta$, it holds that $\chi_{ss}'(G)\le 2\Delta+C$.	
\end{restatable}

Note that no known upper bounds for the strong edge coloring of planar graphs with girth at least some $g$ allows  deriving a bound of the form $2\Delta+C$ for some constant $C$ (see \cite{HLSS2014,BHHV2014,CKKR2018,GZZ2021,CMPR2014,WZ2018,CDYZ2019,BP2013}). Thus, these bounds cannot be used to settle Conjecture \ref{conj-planar-semi-L} even in the case of planar graphs with sufficiently large girth.

 \begin{table}
     \centering
     \begin{tabular}{|m{.27\textwidth}|m{.16\textwidth}|m{.2\textwidth}|m{.25\textwidth}|}
           \hlx{c{2,3,4}v[0]v[0]}
       \multicolumn{1}{c|}{} & \hfill{lower bound}\hfill & \hfill{upper bound} \hfill& \hfill{conjectured value}\hfill\\
           \hlx{v[0]v[0]c{2,3,4}v[0,1,2,3]hvv}
             $\sup\,\{\mathrlap{a'(G)}\phantom{\chi'_{ss}(G)}\colon \Delta(G)=\Delta\}$  & $\phantom{2}\Delta+2$ \hfill\cite{ASZ2001} &  $\phantom{2}\Delta+5$\hfill \cite{SL2024} &   $\phantom{2}\Delta+2$ \hfill\cite{F1978,ASZ2001}\\
               \hlx{vvhvv}
    $\sup\,\{\chi'_{ur}(G)\colon \Delta(G)=\Delta\}$   &  
    { $2\Delta+4$\hfill~\newline (Example~\ref{example})}  &  
    {$2\Delta+4$ if $\g(G)\ge7$  \hfill\newline (Corollary \ref{cor:planar})} &  $2\Delta+4$ \hfill (Conjecture~\ref{our-conj-planar-ur})\\
        \hlx{vvhvv}
    $\sup\,\{\chi'_{ss}(G)\colon \Delta(G)=\Delta\}$   &  
    { $2\Delta+4$\hfill~\newline (Example~\ref{example})}  &  
    {$2\Delta+4$ if $\g(G)\ge7$  \hfill\newline (Theorem \ref{Main-th-girth})} & $2\Delta+C$ \hfill\cite{LMS2024}\newline 
    \\
      \hlx{vvhvv}
         $\sup\,\{\mathrlap{\chi'_s(G)}\phantom{\chi'_{ss}(G)}\colon \Delta(G)=\Delta\}$ &  $4\Delta-3$ \hfill \cite{FSGT1990,WSWC2018}&  $4\Delta+4$ \hfill\cite{FSGT1990}  & \phantom{$2\Delta$}--- \\
           \hlx{vvh}
     \end{tabular}
     \caption{Known, new, and conjectured results on planar graphs.}
     \label{tab:bound}
 \end{table}
 
 \subsection*{Our contribution}
 
\begin{figure}[h!t]
	\begin{center}	\includegraphics[width=.5\textwidth]{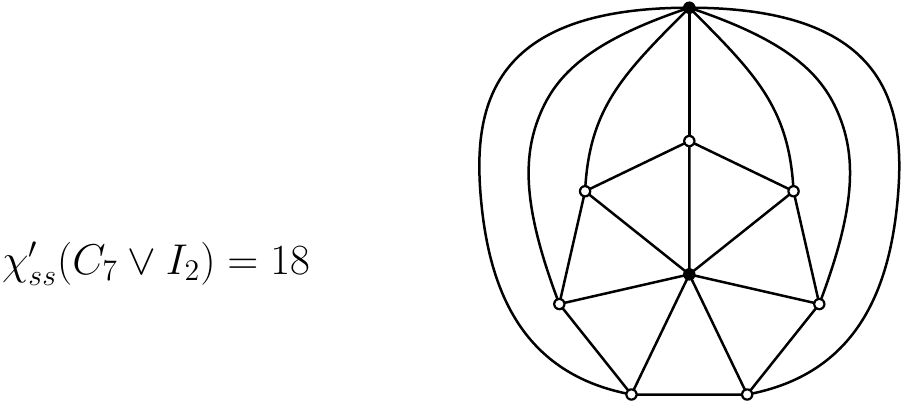}
		\caption{The planar graph $C_7\vee I_2$, which is the join of a $7$-cycle and an independent set of size~$2$, has maximum degree $\Delta=7$ and semistrong chromatic  index $2\Delta+4=18$. Note that every planar graph with maximum degree $\Delta=7$ including $C_7\vee I_2$ as a subgraph has semistrong chromatic  index at least $2\Delta+4$.}
		\label{fig:C7I2}
	\end{center}
\end{figure}

 In this paper, we first provide some lower bounds for the maximum 
 uniquely restricted chromatic index and the maximum semistrong chromatic index of planar graphs with  maximum degree $\Delta\geq 4$ (See \Cref{fig:C7I2} and \Cref{example}), thus showing that  the upper bound for $\chi_{ss}'(G)$ has to be at least of 
the form $2\Delta+C$, for some $C\geq 4$ . 

The main result of this paper bounds the semistrong chromatic index of  graphs with small maximum average degree. Recall that the \emph{maximum average degree} of a graph $G$, denoted $\mad(G)$, is the maximum of the average degrees of the (non-empty) subgraphs of $G$.
Precisely, we prove 
\begin{theorem}\label{Main-th-mad}
	Let  $G$ be a  graph with maximum degree $\Delta$. Then, 
	\begin{align*}
		\chi_{ss}'(G)&\le 2\Delta+4\quad\text{if }\mad(G)<14/5\\
		\text{and }\qquad \chi_{ss}'(G)&\le 2\Delta+2\quad\text{if } \mad(G)<8/3.&
	\end{align*}
\end{theorem}
\begin{corollary}\label{cor:mad}
	Let  $G$ be a  graph with maximum degree $\Delta$. 
	
	Then, 
	$\chi_{ur}'(G)\le 2\Delta+4$ if $\mad(G)<14/5$ and $\chi_{ur}'(G)\le 2\Delta+2$ if $\mad(G)<8/3$.
\end{corollary}

It is an easy consequence of Euler's Formula that a planar graph $G$ with girth $g(G)$ has maximum average degree less than  ${2\g(G)}/({\g(G)-2})$.
Thus, the next result directly follows from Theorem~\ref{Main-th-mad}.

 \begin{theorem}\label{Main-th-girth}
	Let  $G$ be a planar graph with maximum degree $\Delta$. 
	
	Then, $\chi_{ss}'(G)\le 2\Delta+4$ if  $\g(G)=7$ and $\chi_{ss}'(G)\le 2\Delta+2$ if $\g(G)\ge 8$.
\end{theorem}


%
%

\begin{corollary} \label{cor:planar}
    	Let  $G$ be a planar graph with maximum degree $\Delta$. 
        
        Then, $\chi_{ur}'(G)\le 2\Delta+4$ if  $\g(G)=7$ and $\chi_{ur}'(G)\le 2\Delta+2$ if $\g(G)\ge 8$.
\end{corollary}

Note  that  the conjectured upper bound for $a'(G)$ is $\Delta+2$ (the AECC \cite{F1978,ASZ2001}) and the best known upper bound for  general graphs is  $\lceil 3.74(\Delta-1) \rceil+1$  \cite{GKPT2017}.
The two upper bounds in Corollary \ref{cor:mad} are  closer to the AECC for these two classes of sparse graphs with relatively small maximum average degree.

\medskip
\medskip
\noindent{\textbf{Outline.}} 
In \Cref{sec:notation}, we start by introducing some notation and preliminaries. 
 Some example is provided in \Cref{sec:example}, which gives lower bounds for the maximum uniquely restricted chromatic index and semistrong chromatic index of planar graphs with given maximum degree~$\Delta$.
 The overall proof strategy for \Cref{Main-th-mad} is outlined in \Cref{sec:overview}.
 Then, the two cases of 
\Cref{Main-th-mad} are  proved  in \Cref{sec:8/3,sec:14/5}, respectively. 
Finally, we discuss 
some  further research directions in the last section.

 \section{Notation and preliminaries}
 \label{sec:notation}
 Let $G=(V(G),E(G))$ be a  graph.
  For $v\in V(G)$, let $N_{G}(v)=\{u\in V(G): uv\in E(G)\}$ be the open neighborhood of $v$ and  $d_{G}(v)=|N_{G}(v)|$ be the degree of $v$. Let  $\Delta(G)$ 
  denote the maximum degree of $G$ and $\delta(G)$ 
  denote the minimum degree of $G$.  
  The  distance between two edges $e$ and $f$ of $G$ is equal to the distance between their corresponding vertices in the line graph $L(G)$ of $G$, denoted by $d_{G}(e,f)$. 
 
 For a subset $X\subseteq V(G)$, we use $G-X$ to denote the subgraph of $G$ obtained by deleting vertices in $X$. 
 And for a subset $E'\subseteq E(G)$, we use $G\setminus E'$ to denote the subgraph of $G$ obtained by deleting edges in $E'$.
 In particular, we use the abbreviations  $G-v$  and  $G\setminus e$ for the graphs $G-\{v\}$ and $G\setminus \{e\}$, respectively.
 
 Let $i$ be a nonnegative integer.
 An \emph{$i$-vertex}  is a vertex of degree $i$  in $G$.  And  an \emph{$i^{+}$-vertex} (resp.  \emph{$i^{-}$-vertex}) is a vertex of degree at least $i$ (resp. at most $i$)  in $G$.
 For any $uv\in E(G)$ with $d_{G}(u)=i$ (resp. $d_{G}(u)\ge i$ and $d_{G}(u)\le i$), we say that $u$ is an \emph{$i$-neighbor} (resp. \emph{$i^{+}$-neighbor} and \emph{$i^{-}$-neighbor}) of $v$.
 
 For each $v\in V(G)$, let $N_{G}^{1}(v)$ denote the set of 1-neighbors of $v$ in $G$, and let $E_{G}^{1}(v)$ denote the set of edges of $G$ incident  with  $v$ and with the other endvertices in $N_{G}^{1}(v)$.

 We use  $\phi$, $\psi$, $\sigma$ to denote edge colorings. Given   an edge coloring $\phi$ of $G$. For a subset $E'\subseteq E(G)$, we denote by $\phi(E')$ the set of colors that appear on the edges in $E'$. For  two positive integers $i$ and $j$ with $i<j$, we use the abbreviation $[i,j]$ for  the  set of consecutive integers  $\{i,i+1, \dots, j\}$.

 Given a positive integer $k$.
 An edge coloring of $G$ is \emph{good}, if it is a semistrong edge coloring of $G$ using at most $k$ colors. A \emph{good partial coloring} of $G$ is a good coloring $\phi$ of some subgraph $H$ of $G$ such that
 each color class of $\phi$ is a semistrong matching of $G$.  It should be pointed out that 
 in \Cref{sec:8/3}, a good coloring of a graph $G$ is a semistrong edge coloring using at most $2\Delta+2$ colors; 
 and  in  \Cref{sec:14/5}, a good coloring of a graph $G$ is a semistrong edge coloring using at most $2\Delta+4$ colors. 
 
 Let $\phi$ be  a good partial coloring of $G$. Let $e\in E(G)$ be an uncolored edge under $\phi$.
 We say that  a color $\alpha$ is \emph{forbidden} for  $e$ if coloring $e$ with $\alpha$  violates the requirement of semistrong edge coloring.
 And if a color $\alpha$  is not  forbidden for  $e$, then we say that  $\alpha$ is \emph{available} for  $e$ under $\phi$.
 Denote by $F_{\phi}(e)$  and $A_{\phi}(e)$ 
 the sets of colors that are forbidden and available for $e$  under  $\phi$, respectively. 
 It is clear that $A_{\phi}(e)\cup F_{\phi}(e) =[1,k]$ for each uncolored edge $e$ under $\phi$.

 Similarly, if a color $\alpha$ appears on the edges that are at distance at most two from $e$, then we say that $\alpha$  is \emph{strongly forbidden} for  $e$   under $\phi$.
 A color   $\alpha$ is \emph{strongly available} for  $e$ under $\phi$ if $\alpha$  is not strongly forbidden for  $e$.
 We use $\textit{SF}_{\phi}(e)$  and $\textit{SA}_{\phi}(e)$  to denote
 the sets of colors that are  strongly forbidden and strongly  available for $e$  under  $\phi$, respectively. 
 Also, it holds that  $\textit{SA}_{\phi}(e)\cup \textit{SF}_{\phi}(e) =[1,k]$ for each uncolored edge $e$ under $\phi$.
 Moreover, we have $F_{\phi}(e)\subseteq \textit{SF}_{\phi}(e)$ and  $\textit{SA}_{\phi}(e)\subseteq A_{\phi}(e)$ for each uncolored edge $e$ under $\phi$.

\section{Lower bounds}
\label{sec:example}
In this section we provide, through an example, lower bounds for the maximum uniquely restricted chromatic index and semistrong chromatic index of planar graphs with given maximum degree~$\Delta$.

\begin{proposition}\label{example}
	Given an integer $\Delta\ge4$, let $C_{\Delta}\vee I_{2}$ denote the join of the cycle  $C_{\Delta}$ with $\Delta$ vertices and the independent set $I_{2}$ with two vertices. Then,
	\begin{equation*}
		\chi_{ur}'(C_{\Delta}\vee I_{2})=\begin{cases}
			2\Delta+4, \ \text{if} \  \Delta=4,\\
			2\Delta, \ \text{otherwise,}
		\end{cases} 
		\text{and} \ \ 
		\chi_{ss}'(C_{\Delta}\vee I_{2})=\begin{cases}
			2\Delta+4, \ \text{if} \  \Delta=4 \text{ or } \Delta=7,\\
			2\Delta+3, \ \text{otherwise.}
		\end{cases}  
	\end{equation*}
	Moreover, there exist infinitely many planar graph with maximum degree $7$ and semistrong chromatic index at least $18$.
\end{proposition}
\begin{proof}
	Any two edges incident to a vertex in $I_2$ are either adjacent, or belong to a common $C_4$. Hence, these $2\Delta$ edges get different colors in every uniquely restricted coloring and in every semistrong edge coloring. 
   (Refer to \Cref{fig:C7I2} for the case $\Delta=7$.)
	
	Concerning the uniquely restricted coloring of $C_\Delta\vee I_2$ for $\Delta=4$, remark that every edge in $C_\Delta$ belongs to a common $C_4$ with some edge in the complete join. It follows that 
	$\chi_{ur}'(C_4\vee I_2)=2\Delta+\chi_{ur}'(C_4)=2\Delta+4$. However, if $\Delta\geq 5$, then we injectively associate to 
	each edge $e$ of $C_\Delta$  an edge $f$ in the complete join such that $e$ and $f$ do not belong to a $C_4$ and assign to $e$ the color of $f$. Thus, $\chi_{ur}'(C_\Delta\vee I_2)=2\Delta$ in this case.

	As every edge in the complete join as its endpoint in $I_2$ adjacent to the two endvertices of each edge in $C_\Delta$, no edge of $C_\Delta$ can  use a color present in the complete join in a semistrong coloring. Hence, $\chi_{ss}'(C_\Delta\vee I_2)=2\Delta+\chi_{ss}'(C_\Delta)$.
    Note that  it was proved  in \cite{LL2023} that 
    $\chi_{ss}'(C_4)=\chi_{ss}'(C_7)=4$ and $\chi_{ss}'(C_\Delta)=3$ for any other cycle $C_\Delta$. Thus, we obtained the values of $\chi_{ss}'(C_\Delta\vee I_2)$.
	
	Finally, note that every planar graph with maximum degree $7$ including $C_7\vee I_2$ as a subgraph will have semistrong chromatic index at least $\chi_{ss}'(C_7\vee I_2)=18$.
\end{proof}

\section{Proof overview of the main theorem}
\label{sec:overview}

The proof of Theorem~\ref{Main-th-mad} is by contradiction.
In each case, we consider a  counterexample $G$, and the graph $G^{*}$  obtained  from $G$ by  deleting all its vertices of degree 1, and deduce a contradiction by extending a semistrong edge  coloring of a subgraph (or a modified subgraph) of $G$ to the whole graph. 
(Note, instead of introducing $G^*$ we could use the notion of \emph{ghost vertices} introduced in \cite{bonamy2015global}.)
A key property of the considered   (modified) subgraphs  is that  they also  have maximum average degree less than ${14}/{5}$ (resp. ${8}/{3}$) and that any semistrong matching of them has a restriction to $G$ that is  also  a semistrong matching of $G$.
Then, we complete the proof using the discharging method.

The next claim shows that it is sufficient to prove the case $\Delta\ge4$ for Theorem   \ref{Main-th-mad}.

\begin{claim}\label{claim1}
	Theorem~\ref{Main-th-mad} holds for $\Delta\le3$. 
\end{claim}
\begin{clproof}
	Lu{\v{z}}ar, Mockov{\v{c}}iakov{\'a} and Sot{\'a}k  \cite{LMS2024} 
	proved that $\chi_{ss}'(G)\leq8$ for every connected graph $G$ with maximum degree 3 except the complete bipartite graph $K_{3,3}$ and
	$\chi_{ss}'(K_{3,3})=9$, and recently  we proved that  $\chi_{ss}'(C_{4})=\chi_{ss}'(C_{7})=4$ and $\chi_{ss}'(G)\leq3$  for
	any other connected graph with  maximum degree  2 \cite{LL2023}. It follows that Theorem   \ref{Main-th-mad} holds for  the case $\Delta\le3$. 
\end{clproof}

\section{The case of maximum average degree less than ${8}/{3}$} \label{sec:8/3}

In this section,  we prove that if  $G$ is  a  graph with maximum degree $\Delta$ and maximum average degree less than ${8}/{3}$,  then $\chi_{ss}'(G)\le 2\Delta+2$.

Before starting the proof, we  briefly describe the basic idea.
We consider a (potential) counterexample  $G$  with $|V(G)|+|E(G)|$ minimum, that is, $G$ is a connected  graph with maximum degree $\Delta\ge4$ (due to  Claim~\ref{claim1}) and maximum average degree less than ${8}/{3}$ such that $\chi_{ss}'(G)> 2\Delta+2$  and $|V(G)|+|E(G)|$ is as small as possible.
In this section, we say that such a counterexample is \emph{minimal}.

Then, we show the nonexistence of $G$ by showing that the graph $G^*$ obtained from $G$ by deleting all its  1-vertices does not exist.
To do that,  we view each vertex of $G^*$ has an initial charge as its degree,  and then  we show that the special structure of $G^*$ allow charge to be moved from its
 $3^{+}$-vertices  to its $2$-vertices and  “bad” $3$-vertices so that the final charge of each vertex is at least ${8}/{3}$, which leads to a contradiction since $G^*$  has average degree less than ${8}/{3}$.
 
\begin{figure}[htbp]  
	\centering
	\resizebox{8.5cm}{1.8cm}{\includegraphics{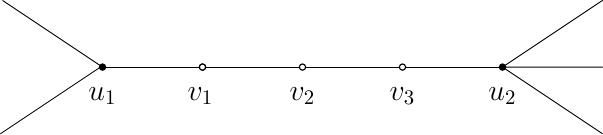}}
	\caption{A $3$-thread.}
		\label{fig:thread}
\end{figure}

We now introduce the concept of a thread in a graph, which plays a crucial role in characterizing the structure of $G^*$.
 Given a positive integer $l$, an \emph{$l$-thread} in a graph $G$ is a trail $u_{1}v_{1}v_{2}\ldots v_{l}u_{2}$ such that $v_{i}$ is a 2-vertex in $G$ for each $i\in[1,l]$ and both  $u_{1}$ and $u_{2}$ are $3^{+}$-vertices in $G$, and we call $u_{1}$ and $u_{2}$  the \emph{ends} of this thread.  See \Cref{fig:thread} for an example.
 If   $u_{1}$ is an end of an  $l$-thread, then we say that \emph{$u_{1}$ is incident with that $l$-thread}. Note that the ends of an $l$-thread with $l\ge2$ may be the same vertex.
 
The structure of $G^*$ is characterized  through a series of lemmas, organized as follows. %
\begin{itemize}
    \item   \Cref{lemma:delta} proves that $G^*$ has no $1$-vertices.
   \item  \Cref{lemma:3+-thread,lemma:2-thread,lemma:1-thread}  describe the properties of $l$-threads in $G^*$,
serving as key tools for proving the subsequent lemmas.
   \item In \Cref{lemma:3-vertex,lemma:3-vertex-bad,lemma:3-vertex-notbad},  different types of $3$-vertices of $G^*$  are analyzed to quantify the charge they may send or receive during discharging.
    \item  \Cref{lemma:4+-vertex} determines  how many vertices a $4^{+}$-vertex   of $G^*$ may be required to help by giving charge.
\end{itemize}

\begin{lemma}\label{lemma:delta}
	$\delta(G^{*})\ge2$.
\end{lemma}
\begin{proof}
	Assume not. Then, 	$\delta(G^{*})\le1$. If 	there is a 0-vertex $v$ in $G^{*}$, then $G$ is isomorphic to the bipartite graph $K_{1,\Delta}$ which has semistrong chromatic index $\Delta$, a contradiction.
	If 	there is a 1-vertex $v$ in $G^{*}$,
	then $v$ has a 1-neighbor $u$ in $G$ as $d_{G}(v)\ge2$. By the minimality of $G$, the graph $G-u$ has a good coloring $\phi$, which is also a good partial coloring of $G$ with the edge $uv$ being uncolored. 	It is clear that  $|\textit{SF}_{\phi}(uv)|\le2\Delta-2$ and thus $|\textit{SA}_{\phi}(uv)|\ge4$. Therefore,  we can color the edge $uv$ properly to extend the coloring $\phi$ to a good  coloring  of $G$, again a contradiction.
\end{proof}

\begin{lemma}\label{lemma:3+-thread}
$G^*$ has no $l$-threads with $l\ge3$.
\end{lemma}
\begin{proof}
Suppose on the contrary that 	$G^{*}$ has an $l$-thread  $u_{1}v_{1}v_{2}\ldots v_{l}u_{2}$ with $l\ge3$. 
	If  $\cup_{i=1}^{l}N_{G}^{1}(v_{i})=\emptyset$,
then each   $v_{i}$ is $2$-vertex in $G$ for $i\in[1,l]$.
		Since $G$ is a minimal counterexample,  the graph $G-v_{2}$ has a good  coloring $\phi$,  which is also a good partial coloring of $G$  with two uncolored edges $v_{1}v_{2}$ and $v_{2}v_{3}$.
		It is easy to check that $|\textit{SF}_{\phi}(v_{1}v_{2})|\le\Delta+1$
	and $|\textit{SF}_{\phi}(v_{2}v_{3})|\le\Delta+1$.
	Therefore, $|\textit{SA}_{\phi}(v_{1}v_{2})|\ge\Delta+1$
	and $|\textit{SA}_{\phi}(v_{2}v_{3})|\ge\Delta+1$.
The coloring $\phi$ can be easily extended to a good coloring   $\psi$ of $G$, a contradiction.	

And if $\cup_{i=1}^{l}N_{G}^{1}(v_{i})\neq\emptyset$,  by the minimality of $G$, the graph $G- \cup_{i=1}^{l}N_{G}^{1}(v_{i})$ has a good coloring $\psi$, which is also a good partial coloring of $G$ with only edges in $\cup_{i=1}^{l}E_{G}^{1}(v_{i})$ being left uncolored.
	It is straightforward to check that $|\textit{SA}_{\psi}(e)|\ge\Delta$ for each $e\in E_{G}^{1}(v_{1})\cup E_{G}^{1}(v_{l})$ and  $|\textit{SA}_{\psi}(e)|\ge2\Delta-2$ for each $e\in E_{G}^{1}(v_{i})$ with $i\in[2,l-1]$.
	Notice that $|E_{G}^{1}(v_{i})|\le\Delta-2$ for each $i\in[1,l]$.	
	Notice also that 
	any two edges $e\in E_{G}^{1}(v_{i})$  and $f\in E_{G}^{1}(v_{j})$ can be colored with a common strongly available color, where 	 $i,j$ are two distinct intergers in $[1,l]$.
	Therefore, it is easy to get a good coloring of $G$ based on $\psi$ by greedily coloring edges in $\cup_{i=1}^{l} E_{G}^{1}(v_{i})$ using their strongly available colors,  again a contradiciton. The lemma holds.
\end{proof}

\begin{lemma}\label{lemma:2-thread}
	Let $u_{1}v_{1}v_{2}u_{2}$ be a  2-thread in	$G^{*}$.
Then, $d_{G}(v_{i})=d_{G^{*}}(v_{i})=2$ for  $i=1,2$, 
 $u_{1}\neq u_{2}$ and
 $d_{G}(u_{i})=d_{G^{*}}(u_{i})\ge4$ for  $i=1,2$.
\end{lemma}
\begin{proof}
	Firstly, we prove that $d_{G}(v_{i})=d_{G^{*}}(v_{i})=2$ for  $i=1,2$. If not, we may assume that $d_{G}(v_{1})\ge3$, then we have $N_{G}^{1}(v_{1})\neq\emptyset$ and $E_{G}^{1}(v_{1})\neq\emptyset$.
	According to the minimality of $G$, the graph $G-N_{G}^{1}(v_{1})$	has a good coloring $\phi$.
	It is easy to see that $|F_{\phi}(e)|\le\Delta+2$ and so $|A_{\phi}(e)|\ge\Delta$ for each $e\in E_{G}^{1}(v_{1})$.
	Because $|E_{G}^{1}(v_{1})|\le \Delta-2$,  a good coloring of $G$ is  obtained
	based on $\phi$ by greedily coloring edges  in  $E_{G}^{1}(v_{1})$.
	Therefore, we obtain a contradiciton.
	
	Secondly, we prove that $u_{1}\neq u_{2}$. If not, suppose $u_{1}=u_{2}$. By the minimality of $G$, the graph $G\setminus v_{1}v_{2}$ has a good coloring $\phi'$. It is easy to see that 
	$|\textit{SF}_{\phi'}(v_{1}v_{2})|\le\Delta$ and so $|\textit{SA}_{\phi'}(v_{1}v_{2})|\ge\Delta+2$.
	 Therefore, the coloring $\phi'$ can be extended to a good coloring of $G$, a contradiction.

	Thirdly, we prove that $d_{G}(u_{i})\ge4$ for  $i=1,2$. If not, w.l.o.g., we may assume that  $d_{G}(u_{1})\le3$. Since $u_{1}$ is the end of the $2$-thread $u_{1}v_{1}v_{2}u_{2}$ in	$G^{*}$, we have $d_{G}(u_{1})=d_{G^{*}}(u_{1})=3$.
	By the minimality of $G$, the graph $G-v_{1}$ has a good coloring $\psi$,  which is also a good partial coloring of $G$  with the two uncolored edges $u_{1}v_{1}$ and $v_{1}v_{2}$. Since   $d_{G}(u_{1})=3$ and $d_{G}(v_{1})=d_{G}(v_{2})=2$, we have 
	 $|\textit{SF}_{\psi}(u_{1}v_{1})|\le2\Delta+1$
	 and $|\textit{SF}_{\psi}(v_{1}v_{2})|\le\Delta+2$.
	It follows that  $|\textit{SA}_{\psi}(u_{1}v_{1})|\ge1$
	 and $|\textit{SA}_{\psi}(v_{1}v_{2})|\ge\Delta$.
	Therefore, we can color $u_{1}v_{1}$ and $v_{1}v_{2}$ properly to obtain a good coloring of $G$ based on $\psi$, a contradiction.
	 
	 Finally, we prove that   $d_{G}(u_{i})=d_{G^{*}}(u_{i})$
	for $i=1,2$. If not, w.l.o.g., we may assume that  $d_{G}(u_{1})>d_{G^{*}}(u_{1})$. Then $u_{1}$ has a 1-neighbor $w_{1}$ in $G$. Due to the minimality of $G$, the graph $G\setminus v_{1}v_{2}$ has a good coloring $\sigma$. We may assume that $\sigma(u_{1}v_{1})\neq \sigma(v_{2}u_{2})$ and thus $\sigma$ is a good partial coloring of $G$ as otherwise we can exchange the colors of $u_{1}v_{1}$ and $u_{1}w_{1}$.
	It is easy to check that $|\textit{SF}_{\sigma}(v_{1}v_{2})|\le2\Delta$ and  so $|\textit{SA}_{\sigma}(v_{1}v_{2})|\ge2$.
	A good coloring of $G$ is easy to get by coloring $v_{1}v_{2}$ properly based on $\sigma$, again a contradiction. The lemma is proved.
\end{proof}

\begin{lemma}\label{lemma:1-thread}
	Let $u_{1}v_{1}u_{2}$ be a 1-thread in	$G^{*}$. If $d_{G^{*}}(u_{1})\le 4$ or 
	$d_{G^{*}}(u_{2})\le 4$, then $d_{G}(v_{1})=d_{G^{*}}(v_{1})=2$.
\end{lemma}
\begin{proof}
	Assume not. Then, $d_{G}(v_{1})>d_{G^{*}}(v_{1})=2$.
	Let   $w$ be a 1-neighbor of 
	 $v_{1}$  in $G$.
	By the minimality of $G$, the graph $G-w$ has a good coloring $\phi$. Since  $d_{G^{*}}(u_{1})\le 4$ or 
	$d_{G^{*}}(u_{2})\le 4$, $|F_{\phi}(v_{1}w)|\le 2\Delta+1$ and thus $|A_{\phi}(v_{1}w)|\ge 1$.  The edge $v_{1}w$ can be colored properly to extend the coloring $\phi$ to a good coloring of $G$, a contradiction. Thus we have $d_{G}(v_{1})=d_{G^{*}}(v_{1})=2$.
\end{proof}

For each positive integer $i$ and each $v\in V(G^{*})$, let $l_{i}(v)$ denote the number  of  $i$-threads in $G^{*}$ that are incident with $v$.

\begin{lemma}\label{lemma:3-vertex}
Let $v$ be a vertex of $G^*$.
If $d_{G^{*}}(v)=3$, then $l_{1}(v)\le2$ and 
 $l_{i}(v)=0$ for each $i\ge2$.
\end{lemma}
\begin{proof}
	It follows from Lemmas \ref{lemma:3+-thread} and \ref{lemma:2-thread} that  $l_{i}(v)=0$ for each $i\ge2$.
Now we prove that $l_{1}(v)\le2$. If not,  $l_{1}(v)=3$. We may assume that $v_{1},v_{2},v_{3}$ are the three 2-neighbors of $v$ in  $G^{*}$.  According to Lemma \ref{lemma:1-thread}, $d_{G}(v_{i})=2$ for each $i\in[1,3]$.
	Since $G$ is a minimal counterexample, the graph $G-v$ has a good coloring $\phi$. It is easy to check that $|\textit{SA}_{\phi}(vv_{i})|\ge\Delta\ge4$ for each $i\in[1,3]$ and  if $d_{G}(v)>3$ then $|\textit{SA}_{\phi}(e)|\ge 2\Delta-1$ for each $e\in E_{G}^{1}(v)$.
	Notice that $|E_{G}^{1}(v)|\le\Delta-3$,
a good coloring of $G$ can be easily obtained based on $\phi$ by first coloring $vv_{1},vv_{2},vv_{3}$ and then
 the  edges in $E_{G}^{1}(v)$ one by one greedily (if $E_{G}^{1}(v)\neq\emptyset$). This gives a contradiciton.
Therefore,    $l_{1}(v)\le2$. The lemma holds.
\end{proof}

	If $d_{G^{*}}(v)=3$ and $l_{1}(v)=2$, then we say that $v$ is a \emph{bad vertex} in $G^{*}$.
If $uv$ is an edge in  $G^{*}$ and $v$  is a bad vertex, then we say that $v$ is a \emph{bad neighbor}  of $u$ in $G^{*}$.
For each $v\in V(G^{*})$, let  $bn(v)$ denote the number of bad neighbors of $v$ in $G^{*}$.

\begin{lemma}\label{lemma:3-vertex-bad}
Let $v$ be a vertex of $G^*$.	If $d_{G^{*}}(v)=3$ and $v$ is bad, then $d_{G}(v)=d_{G^{*}}(v)=3$ and $bn(v)=0$. 
\end{lemma}
\begin{proof}
		Let $v_{1}$ and $v_{2}$ be the two 2-neighbors of $v$ in  $G^{*}$ and  $v_{3}$  the third neighbor of $v$. According to Lemma \ref{lemma:1-thread}, both $v_{1}$ and $v_{2}$ are 2-vertices in  $G$.
	
	First we prove that $d_{G}(v)=d_{G^{*}}(v)=3$. If not, $v$ has a 1-neighbor $u$ in $G$.
		By the minimality of $G$, the graph $G-u$	has a good coloring $\phi$.
		It is easy to check that $|\textit{SF}_{\phi}(uv)|\le2\Delta$ and so  $|\textit{SA}_{\phi}(uv)|\ge2$.  The coloring   $\phi$ can be extended to a good coloring of $G$, a contradiction.

	Then we prove that $bn(v)=0$. If not,
	$v_{3}$ is a bad vertex. By the argument in  the previous paragraph,  we must have $d_{G}(v_{3})=d_{G^{*}}(v_{3})=3$. 
	We may assume that $w_{1}$ and $w_{2}$ are the two 2-neighbors of $v_{3}$ in $G$. By Lemma \ref{lemma:1-thread}, $d_{G}(w_{i})=d_{G^{*}}(w_{i})=2$ for  $i=1,2$. 
	Due to the minimality of $G$, the graph $G-v_{3}$ has a good coloring $\psi$.
	It is easy to check that  $|\textit{SF}_{\psi}(vv_{3})|\le6$ and $|\textit{SF}_{\psi}(v_{3}w_{i})|\le\Delta+3$ for $i=1,2$. Therefore,  $|\textit{SA}_{\psi}(vv_{3})|\ge2\Delta-4\ge4$ and $|\textit{SA}_{\psi}(v_{3}w_{i})|\ge\Delta-1\ge3$ for $i=1,2$. 
	The coloring $\psi$ can be easily extended to a good coloring of $G$ by greedily coloring $vv_{3},v_{3}w_{1},v_{3}w_{2}$, again a contradiction. The lemma is proved.
\end{proof}

\begin{lemma}\label{lemma:3-vertex-notbad}
Let $v$ be a vertex of $G^*$. If $d_{G^{*}}(v)=3$ and  $v$ is not bad,   then $l_{1}(v)+bn(v)\le1$.
\end{lemma}
\begin{proof}	
Denote by   $v_{1},v_{2},v_{3}$ the three neighbors of $v$ in $G^{*}$.
Suppose on the contrary that   $l_{1}(v)+bn(v)\ge2$.
	Since   $v$ is not a bad vertex in $G^{*}$,  $l_{1}(v)\le1$ and thus we must have $bn(v)\ge1$. 	
	We may assume that 
  $v_{1}$ is  a bad neighbor of  $v$  in $G^{*}$. By Lemma \ref{lemma:3-vertex-bad}, $d_{G}(v_{1})=3$. Let $w_{1}$  and $w_{2}$ be the other two neighbors of $v_{1}$.	
According to Lemma \ref{lemma:1-thread}, $d_{G}(w_{i})=d_{G^{*}}(w_{i})=2$ for  $i=1,2$.
 Due to the minimality of $G$, the graph $G-v_{1}$ has a good coloring $\phi$, which is also a good partial coloring of $G$ with three edges $vv_{1},v_{1}w_{1},v_{1}w_{2}$ being uncolored.

Since  $l_{1}(v)+bn(v)\ge2$, at least one of $v_{2}$ and $v_{3}$ is either $2$-vertex or bad vertex in $G^{*}$. We may   assume that $v_{2}$ is such a vertex. If $d_{G^{*}}(v_{2})=2$, then  $l_{1}(v)\ge1$, by Lemma \ref{lemma:1-thread}, $d_{G}(v_{2})=d_{G^{*}}(v_{2})=2$ as  $d_{G^{*}}(v)=3$. And if $v_{2}$ is bad, then by Lemma \ref{lemma:3-vertex-bad}, $d_{G}(v_{2})=d_{G^{*}}(v_{2})=3$. 
In both cases, it is 	easy to check that  $|\textit{SF}_{\phi}(vv_{1})|\le\Delta+5$ and $|\textit{SF}_{\phi}(v_{1}w_{i})|\le\Delta+3$ for $i=1,2$. Thus we have  $|\textit{SA}_{\phi}(vv_{1})|\ge\Delta-3\ge1$ and $|\textit{SA}_{\phi}(v_{1}w_{i})|\ge\Delta-1\ge3$ for $i=1,2$. 
Therefore,	the coloring $\phi$ can be easily extended to a good coloring of $G$ by greedily coloring the three edges $vv_{1},v_{1}w_{1},v_{1}w_{2}$, a contradiction.  The lemma holds.
\end{proof}

\begin{lemma}\label{lemma:4+-vertex}
Let $v$ be a vertex of $G^*$. 	If $d_{G^{*}}(v)\ge4$,   then $$\sum_{i=3}^{d_{G^{*}}(v)}il_{i}(v)+2l_{2}(v)+l_{1}(v)+bn(v)\le2d_{G^{*}}(v)-4.$$
\end{lemma}
\begin{proof}
According to Lemma \ref{lemma:3+-thread}, we have $l_{i}(v)=0$ for each $i\in[3,d_{G^{*}}(v)]$.	
We now  prove that $2l_{2}(v)+l_{1}(v)+bn(v)\le2d_{G^{*}}(v)-4$.
We may assume that $l_{2}(v)>0$ as otherwise the lemma is trivial.
For $i\in[1,l_{2}(v)]$, let $vv_{i}^{1}v_{i}^{2}u_{i}$	denote the $i$-th $2$-thread in $G^{*}$ incident with  $v$ and let
$e_{i}=vv_{i}^{1}$ and $f_{i}=v_{i}^{1}v_{i}^{2}$.	
And if  $l_{2}(v)<d_{G^{*}}(v)$, let $w_{1},w_{2},\ldots,w_{d_{G^{*}}(v)-l_{2}(v)}$ be the other neighbors of $v$ that is not on the $2$-threads.
 Since  $l_{2}(v)>0$, according to Lemma \ref{lemma:2-thread},  $v\neq u_i$, $d_{G}(v)=d_{G^{*}}(v)\ge4$, $d_{G}(u_{i})=d_{G^{*}}(u_{i})\ge4$  and $d_{G}(v_{i}^{j})=d_{G^{*}}(v_{i}^{j})=2$ for each $i\in[1,l_{2}(v)]$ and $j=1,2$.
 It follows that
if $l_{2}(v)\ge2$, then for  any two distinct integers $i,j\in[1,l_{2}(v)]$, 
$f_{i}$ and $f_{j}$ are   at distance 3 in $G$ no matter whether 
 $u_{i}$ is the same as $u_{j}$. 
 
 \begin{claim}\label{claim:4+}
 	$l_{2}(v)+l_{1}(v)\le d_{G}(v)-1$ and $l_{2}(v)\le d_{G}(v)-2$.
 \end{claim}
	\begin{clproof}
		First we prove that  $l_{2}(v)+l_{1}(v)\le d_{G}(v)-1$. If not,  $l_{2}(v)+l_{1}(v)= d_{G}(v)$. That is, each neighbor of $v$ is a $2$-vertex in $G^{*}$.
	By the minimality of $G$, the graph $G-v_{1}^{1}$ has a good coloring $\phi$, which is also a good coloring of $G$ with the two uncolored edges $e_{1}$ and $f_{1}$. It is easy to check that $|F_{\phi}(e_{1})|\le2\Delta-1$ and $|F_{\phi}(f_{1})|\le2\Delta-1$. Therefore,  $|A_{\phi}(e_{1})|\ge3$ and $|A_{\phi}(f_{1})|\ge3$.  Based on $\phi$, a good coloring of $G$ is easy to get by coloring $e_{1}$ and $f_{1}$ properly, a contradiction.	
	
Then we prove that  $l_{2}(v)\le d_{G}(v)-2$. If not,  we must have $l_{2}(v)= d_{G}(v)-1$ as $l_{2}(v)+l_{1}(v)\le d_{G}(v)-1$. Since $G$ is a minimal counterexample, the graph $G-\{v_{1}^{1},v_{2}^{1},\ldots,v_{d_{G}(v)-1}^{1}\}$ has a good coloring $\psi$.
	It is straightforward to check that  $|\textit{SA}_{\psi}(e_{i})|\ge\Delta+1$ and $|\textit{SA}_{\psi}(f_{i})|\ge\Delta+1$ for each $i\in[1,d_{G}(v)-1]$.
	Recall that for any two distinct integers $i,j\in[1,d_{G}(v)-1]$,  $f_{i}$ and $f_{j}$ are at distance 3 in $G$.  We can extend $\psi$ to a good coloring of $G$ by first coloring  $e_{1},e_{2},\ldots,e_{d_{G}(v)-1}$  and then coloring $f_{1},f_{2},\ldots,f_{d_{G}(v)-1}$ one by one greedily. This gives a contradiction.
	\end{clproof}

 Next we prove that $2l_{2}(v)+l_{1}(v)+bn(v)\le 2d_{G}(v)-4$.
If  $l_{2}(v)\le d_{G}(v)-4$, then $2l_{2}(v)+l_{1}(v)+bn(v)\le l_{2}(v)+d_{G}(v)\le 2d_{G}(v)-4$, we are done. Thus 
we may assume that  $d_{G}(v)-3\le l_{2}(v)\le d_{G}(v)-2$ according to 
Claim \ref{claim:4+}.

If  $l_{2}(v)=d_{G}(v)-2$,  then 
$l_{1}(v)+bn(v)\le 2$.
We must have $l_{1}(v)+bn(v)\ge1$ as otherwise $2l_{2}(v)+l_{1}(v)+bn(v)\le 2d_{G}(v)-4$.
It follows that $d_{G^{*}}(w_{1})\le3$ or  $d_{G^{*}}(w_{2})\le3$. 
Since $G$ is a minimal counterexample, the graph $G-\{v_{1}^{1},v_{2}^{1},\ldots,v_{d_{G}(v)-2}^{1}\}$ has a good coloring $\phi$.
It is straightforward to check that 
$|F_{\phi}(e_{i})|\le\Delta+4$ and 
$|\textit{SF}_{\phi}(f_{i})|\le\Delta+2$ for each $i\in[1,d_{G}(v)-2]$.
Therefore, we have
$|A_{\phi}(e_{i})|\ge\Delta-2$ and $|\textit{SA}_{\phi}(f_{i})|\ge\Delta$ for each $i\in[1,d_{G}(v)-2]$.
Recall that for any two distinct integers $i,j\in[1,d_{G}(v)-2]$,  $f_{i}$ and $f_{j}$ are at distance 3 in $G$. 
Based on $\phi$, we can first color $e_{1},e_{2},\ldots,e_{d_{G}(v)-2}$  and then $f_{1},f_{2},\ldots,f_{d_{G}(v)-2}$ one by one greedily. This yields a   good coloring of $G$,  a contradiction.

If $l_{2}(v)=d_{G}(v)-3$, then $l_{1}(v)+bn(v)\le3$. We must have $l_{1}(v)+bn(v)=3$ as otherwise we are done.
It follows that  $d_{G^{*}}(w_{i})\le3$ for each $i\in[1,3]$. 


Now, if $bn(v)=0$,  then  $l_{1}(v)=3$ and so
$w_{1},w_{2},w_{3}$ are all 2-vertices in $G^{*}$.
By the minimality of $G$,  
the graph $G-\{v_{1}^{1},\ldots,v_{d_{G}(v)-3}^{1}\}$ has a good coloring $\psi$. 
It is easy to check that 
$|F_{\psi}(e_{i})|\le7$ and 
$|\textit{SF}_{\psi}(f_{i})|\le\Delta+3$ for each $i\in[1,d_{G}(v)-3]$.
Thus, we have
$|A_{\psi}(e_{i})|\ge2\Delta-5\ge\Delta-1$ and $|\textit{SA}_{\psi}(f_{i})|\ge\Delta-1$ for each $i\in[1,d_{G}(v)-3]$.
Recall that $f_{i}$ and $f_{j}$ are at distance 3 in $G$ for any two distinct integers $i,j\in[1,d_{G}(v)-3]$  if $d_{G}(v)\ge5$. 
We can extend $\psi$ to a good coloring of $G$ by first coloring $e_{1},\ldots,e_{d_{G}(v)-3}$  and then $f_{1},\ldots,f_{d_{G}(v)-3}$ one by one greedily,  a contradiction.

And if $bn(v)>0$,   we may assume that $w_{1}$ is a  bad vertex and $x_{1}$ and $x_{2}$ are the other two neighbors of $w_{1}$. 
By Lemma \ref{lemma:3-vertex-bad}, $d_{G}(w_{1})=d_{G^{*}}(w_{1})=3$.
And according to Lemma \ref{lemma:1-thread}, $d_{G}(x_{i})=d_{G^{*}}(x_{i})=2$ for $i=1,2$.
Due to the minimality of $G$, the graph $G-\{w_{1},v_{1}^{1},\ldots,v_{d_{G}(v)-3}^{1}\}$ has a good coloring $\sigma$.
Recall that  $d_{G^{*}}(w_{2})\le3$ and  $d_{G^{*}}(w_{3})\le3$, 
it is straightforward to check that
$|F_{\sigma}(vw_{1})|\le8$, $|\textit{SF}_{\sigma}(w_{1}x_{1})|\le\Delta+3$, $|\textit{SF}_{\sigma}(w_{1}x_{2})|\le\Delta+3$, 
$|F_{\sigma}(e_{i})|\le7$ and
$|\textit{SF}_{\sigma}(f_{i})|\le\Delta+2$ for each $i\in[1,d_{G}(v)-3]$.
Because $\Delta\ge4$, we have $|A_{\sigma}(vw_{1})|\ge2\Delta-6\ge2$,
$|\textit{SA}_{\sigma}(w_{1}x_{1})|\ge\Delta-1$,
$|\textit{SA}_{\sigma}(w_{1}x_{2})|\ge\Delta-1$,
$|A_{\sigma}(e_{i})|\ge2\Delta-5\ge\Delta-1$ and
$|\textit{SA}_{\sigma}(f_{i})|\ge\Delta$ for each $i\in[1,d_{G}(v)-3]$.
Now we can color $vw_{1}$ and $e_{1},\ldots,e_{d_{G}(v)-3}$ to extend 
 $\sigma$ to a new good partial  coloring $\sigma^{*}$ of $G$, in which $\sigma^{*}(e_{i})\neq \sigma(v_{i}^{2}u_{i})$  for each $i\in[1,d_{G}(v)-3]$.
It is not difficult to see that 
 $|A_{\sigma^{*}}(w_{1}x_{1})|\ge\Delta-2\ge2$,
 $|A_{\sigma^{*}}(w_{1}x_{2})|\ge\Delta-2\ge2$,
$|\textit{SA}_{\sigma^{*}}(f_{i})|\ge2$ for each $i\in[1,d_{G}(v)-3]$.
Notice that for any  $i\in[1,d_{G}(v)-3]$ and $j=1,2$, $f_{i}$ and $w_{1}x_{i}$ are at distance 3 in $G$.
Notice also that any two distinct edges $f_{i}$ and $f_{j}$ are also at distance 3 in $G$ for $i,j\in[1,d_{G}(v)-3]$ if $d_{G}(v)\ge5$.
The coloring $\sigma^{*}$ can be  further extended to a good coloring of $G$ by coloring $w_{1}x_{1},w_{1}x_{2}$
and $f_{1},\ldots,f_{d_{G}(v)-3}$  properly, again a contradiction.
 The lemma is proved.
\end{proof}

We have now all the ingredients to prove the case of maximum average degree less than $8/3$.
\begin{theorem}[Second case of Theorem~\ref{Main-th-mad}]
Let  $G$ be a  graph with maximum degree $\Delta$. 

If $\mad(G)<8/3$, then $\chi_{ss}'(G)\le 2\Delta+2$.
\end{theorem}

\begin{proof}
Assume for contradiction that the Theorem does not holds, and 
let  $G$ be a  counterexample   with $|V(G)|+|E(G)|$ minimum.
Let $G^*$ be the graph obtained from $G$ by deleting all its $1$-vertices.
For each element $v\in V(G^{*})$, we denote by $w(v)$ and $w^{*}(v)$  the  \emph{initial weight} and  the \emph{final weight}, respectively.

Let $w(v)=d_{G^{*}}(v)$ for each $v\in V(G^{*})$.
Because $G$ has maximum average degree less than ${8}/{3}$,  the average degree of  $G^{*}$ does not exceed ${8}/{3}$. 
Hence, we have 
$$\sum_{v\in V(G^{*})}w(v)=\sum_{v\in V(G^{*})}d_{G^{*}}(v)<\frac{8}{3}|V(G^{*})|.$$

The following is the unique discharging rule.

\noindent{\textbf{(R1)}.}  Each $3^{+}$-vertex gives $\frac{1}{3}$ to each 2-vertex on its incident $l$-thread  
and each of  its bad neighbors.

\begin{claim}\label{claim:w*}
	For each $v\in V(G^{*})$, $w^{*}(v)\ge {8}/{3}$.
\end{claim}
\begin{clproof}
By Lemma \ref{lemma:delta}, $G^{*}$ has no 1-vertices.

	If $v$ is a 2-vertex in $G^{*}$, then by \textbf{(R1)}, we  have
	$w^{*}(v)=w(v)+2\times\frac{1}{3}=2+\frac{2}{3}={8}/{3}$.
	
		If $v$ is a 3-vertex in $G^{*}$, then by Lemma \ref{lemma:3-vertex}, it holds that 
		$l_{i}(v)=0$ for each $i\ge2$.
		Now, if   $v$ is bad, by the definition of a bad vertex and Lemma \ref{lemma:3-vertex-bad},  	$l_{1}(v)=2$ and 	$bn(v)=0$.
	By \textbf{(R1)}, we  have
	$w^{*}(v)=w(v)-2\times\frac{1}{3}+\frac{1}{3}=3-\frac{2}{3}+\frac{1}{3}={8}/{3}$.
		And if  $v$ is not bad, by Lemma \ref{lemma:3-vertex-notbad}, we have $l_{1}(v)+bn(v)\le1$. 	By \textbf{(R1)}, we  have
	$w^{*}(v)\ge w(v)-\frac{1}{3}=3-\frac{1}{3}={8}/{3}$.
	
Finally, if $v$ is a $4^{+}$-vertex in $G^{*}$, then by  Lemma  \ref{lemma:4+-vertex}, we have 
$\sum_{i=3}^{d_{G^{*}}(v)}il_{i}(v)+2l_{2}(v)+l_{1}(v)+bn(v)\le2d_{G^{*}}(v)-4$.
	By \textbf{(R1)}, it holds that
$w^{*}(v)
\ge d_{G^{*}}(v)-\frac{1}{3}\times[2d_{G^{*}}(v)-4]\ge \frac{d_{G^{*}}(v)}{3}+\frac{4}{3}\ge{8}/{3}$.

Thus, $w^{*}(v)\ge {8}/{3}$ for each $v\in V(G^{*})$. 
\end{clproof}

It follows from Claim \ref{claim:w*} that, after discharging, we have 
$$\frac{8}{3}|V(G^{*})|\le \sum_{v\in V(G^{*})}w^{*}(v)=\sum_{v\in V(G^{*})}w(v)<\frac{8}{3}|V(G^{*})|,$$
which is a contradiction. 
\end{proof}

 \section{The case of maximum average degree less than ${14}/{5}$} \label{sec:14/5}

 In this section, we prove that if  $G$ is  a  graph with maximum degree $\Delta$ and maximum average degree less than ${14}/{5}$,  then $\chi_{ss}'(G)\le 2\Delta+4$.

 We  consider  a  (potential)  counterexample  $G$ with the minimum number of $2^{+}$-vertices, and subject to this with the minimum number of edges. 
In this section, we say that such a counterexample is \emph{minimal}. (Note that the notion of minimal counterexample used in this section is, in general, different from the one used in the previous section.)
 It is worth mentioning this notion of minimality is inspired by the proofs in \cite{CKKR2018}.

 Similar to  Section \ref{sec:8/3}, 
we will show that the graph $G^*$ obtained from $G$ by deleting all its $1$-vertices does not exist, and thus  $G$ does not exist.
Since  $G$ is a minimal counterexample,
$G$ is a connected  graph with  maximum degree  $\Delta\ge 4$ (due to Claim~\ref{claim1}) and maximum average degree $\mad(G)<{14}/{5}$ such that $\chi_{ss}'(G)> 2\Delta+4$.
 As $G^*$ is actually a subgraph of $G$ induced by all its $2^+$-vertices, 
 $G^{*}$ is  a connected 
graph with  $\mad(G^{*})<{14}/{5}$.
 Denote by $\delta^*$ and $\Delta^{*}$ the minimum and  maximum degrees of the graph  $G^{*}$, respectively. It is obvious that  $\delta^{*}\ge \delta$ and $\Delta^{*}\le \Delta$.

 To show the nonexistence of $G^*$,
 we will first characterize the local properties of
 $G^*$  through a series of lemmas. 
 Then, we will derive a contradiction using degree discharging, where each vertex $v$ of $G^*$ starts with charge $d_{G^*}(v)$,  and the structure of $G^*$ allow the charge 
to be moved  from high-degree vertices to low-degree vertices, such that each vertex ends with a charge of at least ${14}/{5}$.
 
 For readability, we  briefly introduce the structure of the following lemmas.
\begin{itemize}
    \item  \Cref{lem:delta,lem:Delta}  prove that $G^*$ has minimum and maximum degrees at least 2 and 4, respectively.
    \item    \Cref{lem:1-thread} characterizes the edges of $G^*$ whose endpoints both have small degree, and serves as a crucial tool used in nearly every subsequent lemma.

    \item \Cref{lem:3-vertex,lem:3-vertex-bad-neighbor,lem:terrible-3-vertex,lem:terrible-3-vertex-3-neighbor,lem:bad-3-vertex,lem:poor-2-vertex,lem:bad-3-vertex-bad-neighbor}  investigate the properties of  2-vertices and 3-vertices of $G^*$. More specifically, 
    \Cref{lem:poor-2-vertex}  concerns the 2-vertices, while the remaining lemmas focus on the  3-vertices.
    Among them, \Cref{lem:terrible-3-vertex,lem:terrible-3-vertex-3-neighbor,lem:poor-2-vertex,} serve as key tools for proving  subsequent lemmas, while    
    \Cref{lem:3-vertex-bad-neighbor,lem:bad-3-vertex,lem:poor-2-vertex,lem:bad-3-vertex-bad-neighbor} play a particularly crucial role in the discharging process. 
      Note that these lemmas are presented in an order that reflects their logical dependencies, with earlier results serving as tools for later ones.

    \item \Cref{lem:4-vertex,lem:4-vertex-poor-neighbor,lem:Delta-vertex-poor-neighbor,lem:5-vertex,lem:5-vertex-2-poor-neighbor,lem:5-vertex-1-poor-neighbor,lem:6+-vertex-Delta-2,lem:6+-vertex-Delta-3}  analyze how many vertices a $4^+$-vertex in $G^*$ may be required to send charge to during the discharging process. Specifically, \Cref{lem:4-vertex,lem:4-vertex-poor-neighbor} concern the 4-vertices;   \Cref{lem:Delta-vertex-poor-neighbor}  serves as a  tool in subsequent proofs; \Cref{lem:5-vertex,lem:5-vertex-2-poor-neighbor,lem:5-vertex-1-poor-neighbor}  focus on the 5-vertices 
     under the assumption that $\Delta^*\ge5$; and 
     \cref{lem:6+-vertex-Delta-2,lem:6+-vertex-Delta-3} address the $\Delta^*$-vertices  of $G^*$   assuming  $\Delta^*\ge6$.
\end{itemize}

\medskip

 The first lemma is similar to \Cref{lemma:delta}, but the proof slightly differs.
 
 \begin{lemma}\label{lem:delta}
 	$\delta^{*}\ge2$.
 \end{lemma}
 \begin{proof}
 	Assume not. Then, $\delta^{*}\le1$. If 	there is a 0-vertex $v$ in $G^{*}$, then $G$ is isomorphic to the bipartite graph $K_{1,\Delta}$ which has semistrong chromatic index $\Delta$, a contradiction.
 	If 	there is a 1-vertex $v$ in $G^{*}$,
 	then $v$ has a 1-neighbor $u$ in $G$ as $d_{G}(v)\ge2$. Let $H=G-u$. Since $H$ is a subgraph of $G$, $\mad(H)<{14}/{5}$. Moreover,
 the graph	$H$ has no more $2^{+}$-vertices than $G$ and has fewer edges.
 	Thus by the minimality of $G$, the graph $H$ has a good coloring $\phi$, which is also a good partial coloring of $G$ with the edge $uv$ being uncolored. 	It is clear that  $|\textit{SF}_{\phi}(uv)|\le2\Delta-2$ and thus $|\textit{SA}_{\phi}(uv)|\ge6$. Therefore,  we can color the edge $uv$ properly to extend the coloring $\phi$ to a good  coloring $\psi$ of $G$, again a contradiction.
 	The lemma follows.
 \end{proof}

  \begin{lemma}\label{lem:Delta}
 	$\Delta^*\ge4$.
 \end{lemma}
 \begin{proof}
 Suppose on the contrary that	$\Delta^{*}\le3$. Then $G^*$ has a semistrong edge coloring $\phi$ using at most 9 colors  \cite{LMS2024}. 
 According to Lemma \ref{lem:delta}, each vertex of $G^*$ is incident with at most $\Delta-2$ pendant edges in $G$.
 Therefore, we can extend the coloring $\phi$ to a coloring $\psi$ of $G$ by coloring all the pendant edges with the $\Delta-2$ new colors.
Since $\Delta\ge4$, we have $(\Delta-2)+9<2\Delta+4$ and thus
the coloring $\psi$ is a good coloring of $G$. This gives a contradiction.
 \end{proof}

 \begin{lemma}\label{lem:1-thread}
 	Let $uv$ be an edge of $G^{*}$ with $d_{G^{*}}(v)=2$.
 	If $d_{G^{*}}(u)\le 6$,
 then $d_{G}(v)=d_{G^{*}}(v)=2$.
 \end{lemma}
 \begin{proof}
 		If not, $d_{G}(v)>d_{G^{*}}(v)=2$. Then $v$ has a 1-neighbor $w$ in $G$.
 		Consider the graph $H=G-w$. It is clear that   $\mad(H)<{14}/{5}$ and $H$ has no more $2^{+}$-vertices than $G$ and has fewer edges than $G$. 
 	By the minimality of $G$, the graph $H$ has a good coloring $\phi$. Since  $d_{G^{*}}(u)\le 6$, $|F_{\phi}(vw)|\le 2\Delta+3$ and thus $|A_{\phi}(vw)|\ge 1$.  The edge $vw$ can be colored properly to extend the coloring $\phi$ to a good coloring of $G$, a contradiction. The lemma holds.
 \end{proof}

  \begin{lemma}\label{lem:3-vertex}
  Let $v$ be a  $3$-vertex  of $G^*$. Then,  $v$ has at most two $2$-neighhbors in $G^{*}$.
 \end{lemma}
 \begin{proof}
 Suppose on the contraty that $v$ has three $2$-neighbors  $v_{1},v_{2},v_{3}$ in $G^{*}$. 
 By Lemma \ref{lem:1-thread}, $d_{G}(v_{i})=2$ for each $i\in[1,3]$. 
 		Since $G$ is a minimal counterexample, the graph $G-v$ has a good coloring $\phi$. It is easy to check that $|\textit{SA}_{\phi}(vv_{i})|\ge\Delta+2\ge6$ for each $i\in[1,3]$ and  if $d_{G}(v)>3$ then $|\textit{SA}_{\phi}(e)|\ge 2\Delta+1$ for each $e\in E_{G}^{1}(v)$.
 	Notice that $|E_{G}^{1}(v)|\le\Delta-3$,
 	a good coloring of $G$ can be easily obtained based on $\phi$ by first coloring $vv_{1},vv_{2},vv_{3}$ and then
 	the  edges in $E_{G}^{1}(v)$ (if $d_{G}(v)>3$) one by one greedily. This gives a contradiciton.
 	The lemma holds. 
 \end{proof}

To proceed to a more detailed analyses, we classify (thanks to \Cref{lem:3-vertex}) the vertices $v$ of $G^*$ of degree $2$ or $3$ as indicated on \Cref{tab:types}.
\begin{table}[ht]
    \centering
    \renewcommand{\arraystretch}{1.3}
\begin{tabular}{|c|rp{.45\textwidth}|}
\hlx{hvv}
     \multirow{3}*{$3$-vertex}&\emph{good}& $v$ has no $2$-neighbors\\
     &\emph{bad}&$v$ has exactly one $2$-neighbor\\
    &\emph{terrible}&$v$ has exactly two $2$-neighbors\\
    \hlx{vhv}
     \multirow{2}*{$2$-vertex}&\emph{poor}&$v$ has a $2$-neighbor or a terrible $3$-neighbor\\
   &\emph{nonpoor}&$v$ is not poor\\
   \hlx{vvh}
\end{tabular}
    \caption{Different types of $2$- and $3$-vertices of $G^*$.}
    \label{tab:types}
\end{table}

\begin{lemma}\label{lem:3-vertex-bad-neighbor}
  Let $v$ be a  $3$-vertex  of $G^*$. Then,  $v$ has at most two  bad 3-neighbors in $G^{*}$. 
\end{lemma}
\begin{proof}
If not, we assume that $v$ has three   bad 3-neighbors $v_{1},v_{2},v_{3}$.   For each $i\in [1,3]$, let $w_{i}$ denote the 2-neighbor of $v_{i}$ in $G^*$. By Lemma \ref{lem:1-thread}, we have $d_{G}(w_{i})=2$ for each $i\in [1,3]$.
	
	  First we prove that $d_{G}(v_{i})=d_{G^{*}}(v_{i})=3$ for each  $i\in [1,3]$. If not, we may assume that $d_{G}(v_{1})>d_{G^{*}}(v_{1})=3$.  Let $u$ be a
 1-neighbor     of $v_{1}$  in $G$. By the minimality of $G$, the graph $G-u$ has a good coloring $\phi$. It is easy to check that $|F_{\phi}(uv_{1})|\le2\Delta+1$ and so $|A_{\phi}(uv_{1})|\ge3$. We can obtain a good coloring of $G$ based on $\phi$ by coloring $uv_{1}$ properly, a contradiction.

Now by the minimality of $G$, the graph $G-v$ has a good coloring $\psi$. We erase the color of the edge $v_{1}w_{1}$ in $\psi$ to obtain a new good partial coloring $\sigma$ of $G$. Because $d_{G}(v_{i})=d_{G^{*}}(v_{i})=3$ for each  $i\in [1,3]$, it is straightforward to check that $|\textit{SF}_{\sigma}(v_{1}w_{1})|\le 2\Delta$ and $|\textit{SF}_{\sigma}(vv_{i})|\le \Delta+5$ for each  $i\in [1,3]$. And if $d_{G}(v)>d_{G^{*}}(v)$, $|\textit{SF}_{\sigma}(e)|\le 5$  for any edge $e\in E_{G}^{1}(v)$. 
Therefore,  we have $|\textit{SA}_{\sigma}(v_{1}w_{1})|\ge 4$, $|\textit{SA}_{\sigma}(vv_{i})|\ge \Delta-1\ge3$ for each  $i\in [1,3]$ and  $|\textit{SA}_{\sigma}(e)|\ge 2\Delta-1$  for any edge $e\in E_{G}^{1}(v)$  (if $d_{G}(v)>d_{G^{*}}(v)$). 
Notice that $|E_{G}^{1}(v)|\le\Delta-3$, we can extend $\sigma$ to a good coloring of $G$ by first coloring $vv_{1},vv_{2},vv_{3},v_{1}w_{1}$ and then the edges in $E_{G}^{1}(v)$ (if $d_{G}(v)>d_{G^{*}}(v)$). This gives a contradiction.  Thus  $v$ has at most two bad 3-neighbors.
\end{proof}

 \begin{lemma}\label{lem:terrible-3-vertex}
	  Let $v$ be a  $3$-vertex  of $G^*$. 
      If $v$ is terrible, then
 $d_{G}(v)=d_{G^*}(v)=3$.     
      \end{lemma}
\begin{proof}
Suppose on the contrary that $d_{G}(v)>d_{G^{*}}(v)=3$.
Then $v$  has a 1-neighbor $u$ in $G$.
It is obvious that   the graph $G-u$  has maximum average degree less than  ${14}/{5}$ and it has no more $2^{+}$-vertices than $G$ and has fewer edges than $G$.
By the minimality of $G$, the graph $G-u$ has a good coloring $\phi$. 
Since  $v$ is terrible, $v$ has two 2-neighbors $v_1$ and $v_2$ in $G^*$. 
According to  Lemma \ref{lem:1-thread}, we have $d_{G}(v_{1})=d_{G}(v_{2})=2$. 
	It is easy to check that
	$|\textit{SF}_{\phi}(uv)|\le 2\Delta$ and  so $|\textit{SA}_{\phi}(uv)|\ge 4$. We can color the edge $uv$ properly to extend the coloring $\phi$ to a good coloring of $G$, which is  a contradiction. The lemma is proved.
\end{proof}

\begin{lemma}\label{lem:terrible-3-vertex-3-neighbor}
  Let $v$ be a  $3$-vertex  of $G^*$. 
      If $v$ is terrible and  has a  3-neighbor $w$ in $G^{*}$, then $w$ is a good 3-vertex of $G^*$.   
\end{lemma}
\begin{proof}
Suppose on the contrary that $w$ is not good. Then $w$ has a 2-neighbor $w_{1}$.
Since $v$ is terrible, it follows from Lemma \ref{lem:terrible-3-vertex} that $d_{G}(v)=3$.
Denote by $v_{1}$ and $v_{2}$ the other two 2-neighbors of $v$.  By Lemma \ref{lem:1-thread}, $d_{G}(w_{1})=2$ and $d_{G}(v_{1})=d_{G}(v_{2})=2$.
	By the minimality of $G$, the graph $G-v$ has a good coloring $\phi$,  which is also a good partial coloring of $G$ with three uncolored edges $vw,vv_{1},vv_{2}$.	
	 It is easy to check that  	$|\textit{SF}_{\phi}(vw)|\le 2\Delta+1$  and
	 $|\textit{SF}_{\phi}(vv_{i})|\le 2\Delta$ for $i=1,2$.
Therefore, 	$|\textit{SA}_{\phi}(vw)|\ge 3$  and
$|\textit{SA}_{\phi}(vv_{i})|\ge4$ for $i=1,2$.
We can  extend the coloring $\phi$ to a good coloring of $G$ by coloring $vw,vv_{1},vv_{2}$ properly, a contradiction. Therefore, $w$ must be  a good $3$-vertex.		\end{proof}

\begin{figure}[htbp]  
	\centering
	\resizebox{8.1cm}{3.8cm}{\includegraphics{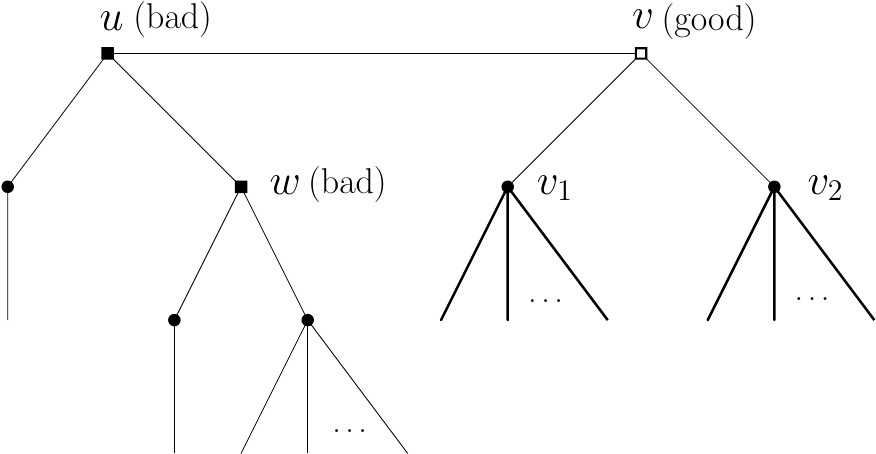}}
	\caption{An illustration of \Cref{lem:bad-3-vertex}, depicting a local structure in $G^*$.}
		\label{fig:good}
\end{figure}

 \begin{lemma}\label{lem:bad-3-vertex}
Let $v$ be a $3$-vertex of $G^{*}$ with three neighbors $u,v_{1},v_{2}$.
If $u$ is a bad 3-vertex  and  it has a bad 3-neighbor $w\notin\{v,v_{1},v_{2}\}$, then  $d_{G^{*}}(v_{i})\ge4$ for  $i=1,2$ (and thus    $v$ is a good 3-vertex). See \Cref{fig:good} for an illustration.
\end{lemma}
\begin{proof}
	If not, we may assume that   $d_{G^{*}}(v_{1})\le3$.
	Let $u_{1}$ be the 2-neighbor of $u$. By Lemma \ref{lem:1-thread}, we have $d_{G}(u_{1})=d_{G^{*}}(u_{1})=2$.

First we prove that $d_{G}(u)=d_{G^{*}}(u)=3$.  If not, $u$ has a 1-neighbor $u'$ in $G$. By the minimality of $G$, the graph $G-u'$ has a good coloring  $\phi$.
Since $u$ is a bad 3-vertex and has a 3-neighbor $v$ in $G^{*}$, 
it  is easy to see that  	$|F_{\phi}(uu')|\le 2\Delta+1$ and so $|A_{\phi}(uu')|\ge 3$. Thus we can extend $\phi$ to a good coloring of $G$ by coloring $uu'$ properly, a contradiction.
Similarly, we also have  $d_{G}(w)=d_{G^{*}}(w)=3$ because $w$  is a bad 3-vertex and has a 3-neighbor $u$.

Then by the minimality of $G$, the graph $G-u$ has a good coloring $\phi$, which is also a good partial coloring of $G$ with three uncolored edges $uv,uu_{1},uw$.
Recall that  $d_{G}(u_{1})=2$, $d_{G}(w)=3$, $d_{G^{*}}(v)=3$ and  $d_{G^{*}}(v_{1})\le3$.  It is straightforward  to check that 
	$|F_{\phi}(uv)|\le 2\Delta+3$, 	$|F_{\phi}(uu_{1})|\le \Delta+4$ and 
	$|F_{\phi}(uw)|\le \Delta+5$.
	Therefore, we have    $|A_{\phi}(uv)|\ge  1$, $|A_{\phi}(uu_{1})|\ge  \Delta\ge4$ and
	  $|A_{\phi}(uw)|\ge \Delta-1\ge3$.
	  By greedily coloring $uv,uw,uu_{1}$ based on $\phi$, we obtain a good coloring of $G$, a contradiction. Therefore,  it holds that $d_{G^{*}}(v_{i})\ge4$ for  $i=1,2$.
\end{proof}


 \begin{lemma}\label{lem:poor-2-vertex}
	Every poor  $2$-vertex of $G^*$ has a $\Delta $-neighbor in $G^{*}$.  This implies that if $G^*$ has a poor $2$-vertex, then $\Delta^*=\Delta\ge4$.
\end{lemma}
\begin{proof}
	Let $u$ be a poor 2-vertex in $G^{*}$. By Lemma \ref{lem:1-thread},  $d_{G}(u)=2$. Let $v$ and $w$ be its two 
	neighbors. 
	 	 We may assume that   $v$ is a 2-neighbor  or a terrible 3-neighbor of $u$.
	 	 	 If $v$ is a  2-neighbor of $u$, by Lemma \ref{lem:1-thread},  we have $d_{G}(v)=d_{G^{*}}(v)=2$ and thus we may assume that  $N_{G}(v)=\{u,v_{1}\}$. 
	 And if $v$ is a terrible 3-neighbor of $u$, 
	by Lemma \ref{lem:terrible-3-vertex}, we must have 	$d_{G}(v)=d_{G^{*}}(v)= 3$ and thus we may assume that
$N_{G}(v)=\{u,v_{1},v_{2}\}$ and $d_{G}(v_{2})=d_{G^{*}}(v_{2})=2$.
Next we prove that $d_{G^{*}}(w)=\Delta$.
	Suppose on the contrary that  
$d_{G^{*}}(w)\le \Delta-1$.
	
We now construst a graph $H$ as follows.	
	If $d_{G}(w)>d_{G^{*}}(w)$, then $w$ has a 1-neighbor $w'$ in $G$,
	  let $H$ be the graph $G\setminus uv$; otherwise,  let  $H$ be the graph obtained from the graph $G\setminus uv$ by  adding	a new pendant edge $ww'$ incident with $w$. 
It is clear that the vertex $w$ always has  two 1-neighbors $u$ and $w'$ in $H$.
As	$d_{G^{*}}(w)\le \Delta-1$, we must have  $\Delta(H)\le\Delta$.
And it is easy to see that
$\mad(H)<{14}/{5}$ and $H$ has fewer $2^{+}$-vertices than $G$.	
By the minimality of $G$, the graph $H$ has a good coloring $\phi$. We may assume that $\phi(uw)\neq\phi(vv_{1})$ as otherwise we can exchange the colors of the two edges $uw$ and $ww'$.	
Now we erase the colors of the edge $ww'$ (if it is not the edge of $G$) and the edge $vv_{2}$ (if $d_{G}(v)=3$) in $\phi$  to get a new good partial coloring $\psi$ of $G$, where the only uncolored edges are $uv$ and $vv_{2}$ (if $d_{G}(v)=3$).

It is easy to check that $|\textit{SF}_{\psi}(uv)|\le 2\Delta+1$
and so  $|\textit{SA}_{\psi}(uv)|\ge 3$. 
And  if $d_{G}(v)=3$, we have $|\textit{SF}_{\psi}(vv_{2})|\le 2\Delta+1$ and so $|\textit{SA}_{\psi}(vv_{2})|\ge 3$. 
Therefore, the coloring $\psi$ can be extended to a good coloring of $G$ easily, a contradiction.
Thus we must have	
$d_{G^{*}}(w)= \Delta$.	
\end{proof}

Recall that each poor 2-vertex of $G^*$
has a 2-neighbor  or a terrible 3-neighbor and  $G^*$ has maximum degree at least $4$ (see Lemma \ref{lem:Delta}),
Lemma \ref{lem:poor-2-vertex} immediately implies the following  easy but useful observation.

 \begin{observation}\label{obser:poor-2-vertex}
 Let $v$ be a $4^{+}$-vertex of $G^*$.
If 	$v$ has a poor 2-neighbor, then  $d_{G^{*}}(v)=\Delta^*=\Delta$.  
\end{observation}

Moreover, it is also easy to observe the following.

\begin{observation}\label{obser:poor-neighbor-nonadjacent}
	If a vertex $u$ of $G^*$ has two poor 2-neighbors  $v$ and $w$, then $vw\notin E(G)$.
\end{observation}
\begin{proof}
	If not,  $vw\in E(G)$. By the minimality of $G$, the graph $G\setminus vw$ has a good coloring $\phi$ with $|\textit{SF}_{\phi}(vw)|\le \Delta$. Thus we can extend $\phi$  to a good coloring of $G$ easily, a contradiction.
\end{proof}

\begin{lemma}\label{lem:bad-3-vertex-bad-neighbor}
 Let $v$ be a  $3$-vertex  of $G^*$. 
      If $v$ is bad, then $v$ 	has at most one bad 3-neighbor and  its 2-neighbor  is not poor,  and it 
	has no terrible 3-neighbors. 
\end{lemma}
\begin{proof}
Let  $v$ be a bad 3-vertex with three neighbors $v_{1},v_{2},v_{3}$ in $G^{*}$. Assume that  $d_{G^{*}}(v_{1})=2$.
 By Lemma \ref{lem:1-thread},  $d_{G}(v_{1})=2$.
	Denote by $v_{1}'$ the other neighbor of $v_{1}$.
	If both $v_{2}$ and $v_{3}$ are bad 3-neighbors, there is a contradiction to Lemma  \ref{lem:bad-3-vertex}. Thus $v$ has at most one bad 3-neighbor.

	If $v_{1}$ is poor, then by the definition of the poor 2-vertex, $v_{1}'$ must be a 2-neighbor or a terrible 3-neighbor of $v_{1}$ (notice that $v$ is a bad 3-vertex).
	But then $v_{1}$ has no $\Delta$-neighbor (notice that $\Delta\ge4$) in $G^{*}$,
		 which is a contradiciton to Lemma \ref{lem:poor-2-vertex}.
	Therefore, $v_{1}$ is not poor. 
	
	Next, we prove that $v$ has no terrible 3-neighbors. If not, we  may assume that $v_{2}$ is a terrible 3-vertex with two 2-neighbors $x$ and $y$. Due to Lemma \ref{lem:terrible-3-vertex}, $d_{G}(v_{2})=3$.
	And by Lemma \ref{lem:1-thread}, $d_{G}(x)=d_{G}(y)=2$.
	By the minimality of $G$, the graph $G-v_{2}$ has a good coloring $\phi$, which is also a good partial coloring of $G$ with three uncolored edges $vv_{2},v_{2}x$ and $v_{2}y$.
	It is easy to check that $|\textit{SF}_{\phi}(vv_{2})|\le 2\Delta+1$, $|\textit{SF}_{\phi}(v_{2}x)|\le 2\Delta$ and  $|\textit{SF}_{\phi}(v_{2}y)|\le 2\Delta$.
	Hence, we have $|\textit{SA}_{\phi}(vv_{2})|\ge 3$, $|\textit{SA}_{\phi}(v_{2}x)|\ge 4$ and  $|\textit{SA}_{\phi}(v_{2}y)|\ge 4$.
	We can extend $\phi$ to a good coloring of $G$ by coloring $vv_{2},v_{2}x,v_{2}y$ properly, a contradiction.
	The lemma is proved.
\end{proof}

Recall that $\Delta^*\ge4$ (see Lemma \ref{lem:Delta}),
 the following lemmas  consider the properties of $4^+$-vertices of $G^*$.

\begin{lemma}\label{lem:4-vertex}
 Let $v$ be  a  $4$-vertex  of $G^*$. 	Then, $v$ has at most three neighbors  that are 2-vertices or bad 3-vertices.
\end{lemma}
\begin{proof}
Let $v_{1},v_{2},v_{3},v_{4}$ be the four neighbors of $v$ in $G^{*}$.  Denote by $e_{i}$ the edge $vv_{i}$ for each $i\in[1,4]$.
Suppose on the contrary that  $v_{1},v_{2},v_{3},v_{4}$  are all  2-vertices or bad 3-vertices.

First we prove that $d_{G}(v_{i})=d_{G^{*}}(v_{i})$ for each $i\in[1,4]$. If $d_{G^{*}}(v_{i})=2$, then by Lemma \ref{lem:1-thread}, the statement is true. And if $d_{G}(v_{i})>d_{G^{*}}(v_{i})=3$, $v_{i}$ has a 1-neighbor $u_{i}$ in $G$. By the minimality of $G$, the graph $G-u_{i}$ has a good coloring $\phi$. Since $v_{i}$ is a bad $3$-vertex, we have $|F_{\phi}(u_{i}v_{i})|\le2\Delta+2$ and  $|A_{\phi}(u_{i}v_{i})|\ge2$. The coloring $\phi$ can be extended to a good coloring of $G$ by coloring $u_{i}v_{i}$ properly, a contradiction. Therefore, in case $d_{G^{*}}(v_{i})=3$,  we also have $d_{G}(v_{i})=d_{G^{*}}(v_{i})=3$.

Then   we prove that $d_{G}(v)=d_{G^{*}}(v)=4$. If not, $v$ has a 1-neighbor $u$ in $G$. By the minimality of $G$, the graph $G-u$ has a good coloring $\psi$. It is easy to see that  $|\textit{SF}_{\psi}(uv)|\le\Delta+7$ and thus
$|\textit{SA}_{\psi}(uv)|\ge\Delta-3\ge1$. Thus we can extend $\psi$ to a good coloring of $G$ by coloring $uv$ properly, a contradiction.


By  the minimality of $G$, the graph $G-v$ has a good coloring $\sigma$.
The remaining proof are divided into the following two cases.
 
\textbf{Case 1.} $E(G[\{v_{1},v_{2},v_{3},v_{4}\}])\neq \emptyset$.

W.l.o.g., we may assume that $v_{1}v_{2}\in E(G)$. If $d_{G}(v_{1})=2$,   then it is easy to check that $|\textit{SF}_{\sigma}(e_{1})|\le6$,   $|\textit{SF}_{\sigma}(e_{2})|\le\Delta+5$ and  $|\textit{SF}_{\sigma}(e_{i})|\le\Delta+6$ for $i=3,4$. 
Thus we have $|\textit{SA}_{\sigma}(e_{1})|\ge2\Delta-2\ge6$,
$|\textit{SA}_{\sigma}(e_{2})|\ge\Delta-1\ge3$,
and  $|\textit{SA}_{\sigma}(e_{i})|\ge\Delta-2\ge2$ for $i=3,4$. 
The coloring $\sigma$ can be extended to a good coloring of $G$ by greedily coloring $e_{4},e_{3},e_{2},e_{1}$ one by one, a contradiction.
Thus $d_{G}(v_{1})=3$. Similarly, we have $d_{G}(v_{2})=3$. Thus both $v_{1}$ and $v_{2}$ are bad 3-vertices.

 Denote by $w_{1}$ and $w_{2}$ the 2-neighbors of $v_{1}$ and $v_{2}$ in $G^{*}$, respectively. By Lemma \ref{lem:1-thread},  $d_{G}(w_{1})=d_{G}(w_{2})=2$.
At this time, we erase the color of $v_{1}v_{2}$ in $\sigma$ yielding a new good partial coloring. This coloring is still denoted by $\sigma$.
 It is not difficult to check that $|\textit{SF}_{\sigma}(v_{1}v_{2})|\le4$, $|\textit{SF}_{\sigma}(e_{i})|\le7$ for $i=1,2$ and $|\textit{SF}_{\sigma}(e_{j})|\le\Delta+6$ for $j=3,4$. 
 Therefore, we have 
$|\textit{SA}_{\sigma}(v_{1}v_{2})|\ge2\Delta\ge8$, $|\textit{SA}_{\sigma}(e_{i})|\ge2\Delta-3\ge5$ for $i=1,2$ and $|\textit{SA}_{\sigma}(e_{j})|\ge\Delta-2\ge2$ for $j=3,4$. 
We can extend $\sigma$ to a good coloring of $G$ by greedily coloring $e_{3},e_{4},e_{1},e_{2},v_{1}v_{2}$ one by one, again a contradiction.

\textbf{Case 2.} $E(G[\{v_{1},v_{2},v_{3},v_{4}\}])= \emptyset$.

In this case,  $v_{i}v_{j}\notin E(G)$ for any $i,j\in[1,4]$.
Now  if  $v_{1},v_{2},v_{3},v_{4}$  are all  2-vertices, it is easy to check that $|\textit{SF}_{\sigma}(e_{i})|\le\Delta+3$ and so $|\textit{SA}_{\sigma}(e_{i})|\ge\Delta+1\ge5$ for each $i\in[1,4]$.
The coloring $\sigma$ can be extended to a good coloring of $G$ by  coloring $e_{1},e_{2},e_{3},e_{4}$ greedily, a contradiciton.

Thus we may assume that  $v_{1},\ldots,v_{p}$ are  bad 3-vertices, where $1\le p\le 4$.
 For each $i\in[1,p]$, let  $w_{i}$ be the 2-neighbor of $v_{i}$ and let $f_{i}$ denote the edge $v_{i}w_{i}$. By Lemma \ref{lem:1-thread}, $d_{G}(w_{i})=2$ for   each $i\in[1,p]$.
If $p\ge2$,  we must have $w_{i}w_{j}\notin E(G)$ for any two integers $i,j\in[1,p]$ as otherwise there is a contradiction to Lemma \ref{lem:poor-2-vertex}.
However,
it is possible that $w_{i}=w_{j}$ for some two integers $i,j\in[1,p]$, but it does not affect the following proof.

We   erase the colors of $f_{1},\ldots,f_{p}$ in $\sigma$ to obtain a new good partial coloring $\sigma'$  of $G$.
It is straightforward to check that
$|\textit{SF}_{\sigma'}(e_{i})|\le\Delta+4$ for each $i\in[1,4]$ and $|\textit{SF}_{\sigma'}(f_{i})|\le2\Delta$ for each $i\in[1,p]$.
 Thus we have  $|\textit{SA}_{\sigma'}(e_{i})|\ge\Delta\ge4$ for each $i\in[1,4]$ and  $|\textit{SA}_{\sigma'}(f_{i})|\ge4$ for each $i\in[1,p]$. We can further extend $\sigma'$ to a new good partial coloring $\sigma^{*}$ of $G$ by greedily coloring $e_{1},e_{2},e_{3},e_{4}$ one by one with their strongly available colors.
 
 Notice that for each $i\in[1,4]$,  the color of the edge  $e_{i}$ is different  from  the colors of   all colored edges that are at distance at most two from it under $\sigma^{*}$.
 Notice also that  $d_{G}(w_{i})=2$
 for each $i\in[1,p]$.
If  $p\ge2$ and $w_{i}=w_{j}$ for some two integers $i,j\in[1,p]$, then both  $f_{i}$ and $f_{j}$ can not be colored with the same color as the  two edges $e_{i}$ and $e_{j}$.
And in other cases,  $f_{i}$  can not be colored with the same color as $e_{i}$ for each $i\in[1,p]$.
Hence, we have $|A_{\sigma^{*}}(f_{i})|\ge2$ for each $i\in[1,p]$.
Recall that  $d_{G}(w_{i})=2$ for each   $i\in[1,p]$,
 there is at most one other edge $f_{j}$ that are at distance at most 2 from $f_{i}$. 
Therefore, we can further extend $\sigma^{*}$ to a good coloring of $G$ by coloring $f_{1},\ldots,f_{p}$ greedily, a contradiction. The lemma is proved.
\end{proof}

 
 \begin{lemma}\label{lem:4-vertex-poor-neighbor}
  Let $v$ be  a  $4$-vertex  of $G^*$. 	Then, 	
  $v$ has at most one poor  2-neighbor. In particular, if  $v$ has exactly one poor  2-neighbor, then $\Delta=4$ and $v$ has at least two $4$-neighbors in $G^{*}$.
 \end{lemma}
\begin{proof}
First  we prove that  $v$ has at most one poor  2-neighbor. If not, let   $v_{1}$ and $v_{2}$ be 
	the two poor  2-neighbors of $v$.  By Observations \ref{obser:poor-2-vertex} and \ref{obser:poor-neighbor-nonadjacent}, we have $\Delta=4$ and $v_{1}v_{2}\notin E(G)$. 
	For $i=1,2$, denote by $w_{i}$ the other neighbor of $v_{i}$. It is clear that  $w_{1}$ and $w_{2}$ are 2-vertices or terrible 3-vertices. It is possible that $w_{1}=w_{2}$, but  it doesn't affect the following proof.
	 By Lemmas \ref{lem:1-thread} and \ref{lem:terrible-3-vertex}, we must have $d_{G}(w_{i})=d_{G^{*}}(w_{i})\le3$. 
	By the minimality of $G$, the graph $G-\{v_{1},v_{2}\}$ has a good colorinng $\phi$.
	Since $d_{G^{*}}(v)=\Delta=4$, it is easy to see that 
	 $|\textit{SF}_{\phi}(vv_{i})|\le2\Delta+2$ and 
	  $|\textit{SF}_{\phi}(v_{i}w_{i})|\le\Delta+4$ for $i=1,2$.
	  It follows that	 $|\textit{SA}_{\phi}(vv_{i})|\ge2$ and 
	  $|\textit{SA}_{\phi}(v_{i}w_{i})|\ge\Delta=4$ for  $i=1,2$.  Therefore, 
coloring $vv_{1},vv_{2},v_{1}w_{1},v_{2}w_{2}$ greedily based on $\phi$ gives rise to a good coloring of $G$, a contradiction.

Suppose that $v$ has exactly one poor 2-neighbor $v_{1}$ and $w_{1}$ is the other neighbor of $v_{1}$.
By the definition of the poor 2-vertex, $w_{1}$ is either a 2-vertex or a terrible 3-vertex. By Lemmas \ref{lem:1-thread} and \ref{lem:terrible-3-vertex}, we  have $d_{G}(w_{1})=d_{G^{*}}(w_{1})\le3$. 
 By Observation \ref{obser:poor-2-vertex},  we have $\Delta=4$.
Next we prove that    $v$ has at least two $4$-neighbors in $G^{*}$. If not, let $v_{2}$ and $v_{3}$ be the other two $3^{-}$-neighbors of $v$ in $G^{*}$. By the minimality of $G$, the graph $G-v_{1}$ has a good coloring $\psi$, which is also a good partial coloring of $G$ with two uncolored edges $vv_{1}$ and $v_{1}w_{1}$.

Now if $w_{1}$ is a 2-vertex, then  $|F_{\psi}(vv_{1})|\le\Delta+7$ and  $|\textit{SF}_{\psi}(v_{1}w_{1})|\le\Delta+3$. 
It follows that  $|A_{\psi}(vv_{1})|\ge\Delta-3=1$ and $|\textit{SA}_{\psi}(v_{1}w_{1})|\ge\Delta+1$. We can color $vv_{1}$ and $v_{1}w_{1}$ properly based on $\psi$ to obtain a good coloring of $G$, a contradiction.

And if $w_{1}$ is a  terrible $3$-vertex, let $x_{1}$ be another 2-neighbor of $w_{1}$. We erase the color of $w_{1}x_{1}$ in $\psi$ to obtain a new good partial coloring $\sigma$ of $G$. It is easy to check that $|F_{\sigma}(vv_{1})|\le\Delta+7$, $|\textit{SF}_{\sigma}(v_{1}w_{1})|\le\Delta+4$ and $|\textit{SF}_{\sigma}(w_{1}x_{1})|\le2\Delta$.  Therefore, we have $|A_{\sigma}(vv_{1})|\ge\Delta-3=1$, $|\textit{SA}_{\sigma}(v_{1}w_{1})|\ge\Delta=4$ and $|\textit{SA}_{\sigma}(w_{1}x_{1})|\ge4$.
Therefore, we can   extend $\sigma$ to a good coloring of $G$ by coloring $vv_{1},v_{1}w_{1},w_{1}x_{1}$ greedily, again a contradiction.  
The lemma holds.
\end{proof}

 \begin{lemma}\label{lem:Delta-vertex-poor-neighbor}
 Let $v$ be  a  $\Delta^*$-vertex  of $G^*$. 	Then,  $v$ has at most $\Delta^*-2$ poor  2-neighbors. 
\end{lemma}
\begin{proof}
	If not, $v$ has at least $\Delta^*-1$ poor  2-neighbors. 
  By observation  \ref{obser:poor-2-vertex}, we have $\Delta^*=\Delta$.
    Let $v_{1},v_{2},\ldots,v_{\Delta-1}$ be the  poor  2-neighbors of $v$. 
	By Observation \ref{obser:poor-neighbor-nonadjacent},
 $v_{i}v_{j}\notin E(G)$ 
for any two distinct integers $i,j\in[1,\Delta-1]$. 
	For each $i\in[1,\Delta-1]$, let $w_{i}$ denote the other neighbor of $v_{i}$ and let $e_{i}=vv_{i}$ and $f_{i}=v_{i}w_{i}$.
	 By the definiton of the poor 2-vertex, each $w_{i}$ is either a 2-vertex or a terrible 3-vertex.  According to Lemmas \ref{lem:1-thread} and \ref{lem:terrible-3-vertex}, we  have $d_{G}(w_{i})=d_{G^{*}}(w_{i})\le3$ for each $i\in[1,\Delta-1]$. 
	 According to Lemmas \ref{lem:terrible-3-vertex-3-neighbor} and \ref{lem:poor-2-vertex}, we must have   $w_{i}w_{j}\notin E(G)$ for any two distinct integers $i,j\in[1,\Delta-1]$.
	 However, it is possible that  $w_{i}= w_{j}$  for some integers $i,j\in[1,\Delta-1]$, but it does not matter.
	 
	 By the minimality of $G$, the graph $G-\{v_{1},v_{2},\ldots,v_{\Delta-1}\}$ has a good coloring $\phi$. It is easy to check that   $|\textit{SF}_{\phi}(e_{i})|\le\Delta+2$ and 
	 $|\textit{SF}_{\phi}(f_{i})|\le\Delta+3$ for each $i\in[1,\Delta-1]$. 
	 Therefore,  $|\textit{SA}_{\phi}(e_{i})|\ge\Delta+2$ and 
	 $|\textit{SA}_{\phi}(f_{i})|\ge\Delta+1$ for each $i\in[1,\Delta-1]$. 
	 We can color $e_{1},e_{2},\ldots,e_{\Delta-1}$ greedily to extend $\phi$ to another good partial coloring $\psi$ of $G$.
		 It is easy to see that  $|\textit{SA}_{\psi}(f_{i})|\ge2$ for each $i\in[1,\Delta-1]$. 
		 Notice that for each  $i\in[1,\Delta-1]$, there is at most one edge $f_{j}$ that is at distance 1 from $f_{i}$ in $G$, where  $j\in[1,\Delta-1]$. Notice  also that any two distinct edges $f_{i}$ and $f_{j}$ that are at distince at least two can be colored with a common strongly available color. We can further extend $\psi$ to a good coloring of $G$ by greedily coloring $f_{1},f_{2},\ldots,f_{\Delta-1}$ one by one, a contradiction. The lemma holds.
\end{proof}

\begin{lemma}\label{lem:5-vertex}
Suppose that $\Delta^*\ge5$.
Let $v$ be  a  $5$-vertex  of $G^*$. 	Then, $v$ has at most two poor  2-neighbors.  
\end{lemma}
\begin{proof}
	 If not,  $v$	 has at least three    poor  2-neighbors.   By Observation \ref{obser:poor-2-vertex}, we must have $\Delta=5$.
	 Then	  according to Lemma \ref{lem:Delta-vertex-poor-neighbor}, $v$	 has exactly three    poor  2-neighbors.	 
	  We may assume that $v_{1},v_{2},v_{3}$ are the three poor  2-neighbors of $v$.   By Observation \ref{obser:poor-neighbor-nonadjacent}, $v_{i}v_{j}\notin E(G)$ for any two distinct integers $i,j\in[1,3]$.
	 For each $i\in[1,3]$, let $w_{i}$ denote the other neighbor of $v_{i}$ and let $e_{i}=vv_{i}$ and $f_{i}=v_{i}w_{i}$.
	 By the definition of the poor 2-vertex, each $w_{i}$ is a 2-vertex or a terrible 3-vertex, where $i\in[1,3]$.
	 By the minimality of $G$, the graph $G-\{v_{1},v_{2},v_{3}\}$ has a good coloring $\phi$, which is also a good partial coloring of $G$ with the six uncolored edges $e_{1},e_{2},e_{3},f_{1},f_{2},f_{3}$.

	 If $w_{1},w_{2},w_{3}$ are all   2-vertices,
	 	  then  by Lemma \ref{lem:1-thread}, 	$d_{G}(w_{i})=2$ for each $i\in[1,3]$. 	  
	 	  It is easy to see that $|\textit{SF}_{\phi}(e_{i})|\le2\Delta+1$ and 
	 $|\textit{SF}_{\phi}(f_{i})|\le\Delta+2$ for each $i\in[1,3]$. 
	 Therefore,  $|\textit{SA}_{\phi}(e_{i})|\ge3$ and 
	 $|\textit{SA}_{\phi}(f_{i})|\ge\Delta+2=7$ for each $i\in[1,3]$.  We can extend $\phi$ to a good coloring of $G$ by greedily coloring the six edges $e_{1},e_{2},e_{3},f_{1},f_{2},f_{3}$ one by one, a contradiction.
	 
	 If exactly one vertex in $\{w_{1},w_{2},w_{3}\}$ is a  terrible 3-vetrex, we may assume that  $w_{1}$ is a  terrible 3-vertex, then it is easy  to check that  $|\textit{SF}_{\phi}(e_{1})|\le2\Delta+2$,  $|\textit{SF}_{\phi}(f_{1})|\le\Delta+4$, $|\textit{SF}_{\phi}(e_{i})|\le2\Delta+1$ and 
	 $|\textit{SF}_{\phi}(f_{i})|\le\Delta+2$ for  $i=2,3$.
	 Thus we have   $|\textit{SA}_{\phi}(e_{1})|\ge2$,  $|\textit{SA}_{\phi}(f_{1})|\ge\Delta=5$, $|\textit{SF}_{\phi}(e_{i})|\ge3$ and 
	 $|\textit{SF}_{\phi}(f_{i})|\ge\Delta+2=7$ for  $i=2,3$.
	 The coloring $\phi$ can be extended to a good coloring of $G$ by greedily coloring the six edges $e_{1},e_{2},e_{3},f_{1},f_{2},f_{3}$, a contradiction.

	Next we  assume that $w_{1},\ldots,w_{p}$ are  terrible 3-vertices, where $2\le p\le 3$. 
	 If there are two distinct integers $i,j\in[1,p]$ such that $w_{i}=w_{j}$, we may assume that $w_{1}=w_{2}$, then it is easy to check that  $|\textit{SF}_{\phi}(e_{3})|\le2\Delta+2$,  $|\textit{SF}_{\phi}(f_{3})|\le\Delta+4$, $|\textit{SF}_{\phi}(e_{i})|\le2\Delta+1$ and 
	 $|\textit{SF}_{\phi}(f_{i})|\le\Delta+2$ for  $i=1,2$.
	 Then we have   $|\textit{SA}_{\phi}(e_{3})|\ge2$,  $|\textit{SA}_{\phi}(f_{3})|\ge\Delta=5$, $|\textit{SF}_{\phi}(e_{i})|\ge3$ and 
	 $|\textit{SF}_{\phi}(f_{i})|\ge\Delta+2=7$ for  $i=1,2$.
	   The coloring $\phi$ can be extended to a good coloring of $G$ by greedily coloring the six edges $e_{3},e_{1},e_{2},f_{3},f_{1},f_{2}$, a contradiction.

	  Thus we now assume that $w_{i}\neq w_{j}$ for   any two distinct integers $i,j\in[1,p]$, where $2\le p\le3$.	 
	  Since  $w_{1},\ldots,w_{p}$ are all  terrible 3-vertices,  by Lemma \ref{lem:terrible-3-vertex-3-neighbor}, $w_{i}w_{j}\notin E(G)$  for   any two distinct integers $i,j\in[1,p]$.
	  This implies that any two edges $f_{i}$
 and $f_{j}$ are at distance 3 in $G$ for any two distinct integers $i,j\in[1,p]$.
 
 	 For each $i\in[1,p]$, let  $x_{i}$ be the other 2-neighbor of $w_{i}$ and let $g_{i}$ denote the edge $w_{i}x_{i}$. 	 	 
	 Since  $w_{1},\ldots,w_{p}$ are  terrible 3-vertices, $x_{1},\ldots,x_{p}$ are all poor 2-vertices and thus $d_{G}(x_{i})=2$ for each $i\in[1,p]$.
	 Moreover,  we must have  $x_{i}\neq x_{j}$ and $x_{i}x_{j}\notin E(G)$ for any two distinct integers $i,j\in[1,p]$ as otherwise there is a contradiction to Lemma \ref{lem:poor-2-vertex}.
	 Therefore, any two edges $g_{i}$ and $g_{j}$ are at distance at least three in $G$, where  $i,j\in[1,p]$.
	
	Now we erase the colors of $g_{1},\ldots,g_{p}$ in $\phi$ yielding a new good partial coloring $\psi$  of $G$.
	  It is straightforward to check that
	 $|\textit{SF}_{\psi}(e_{i})|\le2\Delta+1$,
	 $|\textit{SF}_{\psi}(f_{i})|\le\Delta+3 $
	  for each $i\in[1,3]$ and  $|\textit{SF}_{\psi}(g_{i})|\le2\Delta$ for each $i\in[1,p]$.
	Thus we have
	 $|\textit{SA}_{\psi}(e_{i})|\ge3$,
$|\textit{SA}_{\psi}(f_{i})|\ge\Delta+1=6 $ 
	for each $i\in[1,3]$ and  $|\textit{SF}_{\psi}(g_{i})|\ge4$ for each $i\in[1,p]$.
	We can extend $\psi$ to a new good partial coloring $\sigma$ of $G$ by greedily coloring $e_{1},e_{2},e_{3}$ one by one.
	It is easy to see that   $|\textit{SA}_{\sigma}(f_{i})|\ge\Delta-2=3$
	for each $i\in[1,3]$ and  $|\textit{SF}_{\sigma}(g_{i})|\ge3$ for each $i\in[1,p]$. 
	Then we can  color  $g_{1},\ldots,g_{p}$ greedily to extend $\sigma$ to   a new good partial coloring $\sigma^{*}$  of $G$, in which $|\textit{SA}_{\sigma^{*}}(f_{i})|\ge\Delta-3=2$
	for each $i\in[1,3]$.
	Recall that for any two distinct integers $i,j\in[1,p]$, the two edges $f_{i}$
	and $f_{j}$ are at distance 3 in $G$ and thus they can be colored with a common strongly  available color.
	We can further extend $\sigma^{*}$ to a good coloring of $G$, again a contradiction. The lemma is proved.
\end{proof}

\begin{lemma}\label{lem:5-vertex-2-poor-neighbor}
	Suppose that $\Delta^*\ge5$.
Let $v$ be  a  $5$-vertex  of $G^*$. 
If  $v$ has exactly two poor  2-neighbors, then $v$ has at least two $4^{+}$-neighbors in $G^{*}$.  
\end{lemma}
\begin{proof}
Let $v_{1}$ and $v_{2}$ be the two poor  2-neighbors of $v$. By Lemma \ref{lem:1-thread}, $d_{G}(v_{1})=d_{G}(v_{2})=2$.
 By Observation \ref{obser:poor-2-vertex}, $\Delta=5$. And by Observation \ref{obser:poor-neighbor-nonadjacent}, $v_{1}v_{2}\notin E(G)$. Let $w_{1}$ and $w_{2}$ denote the other neighbors of $v_{1}$ and $v_{2}$, respectively.  By the definition  of the poor 2-vertex, $w_{1}$ and $w_{2}$ are either 2-vertices or terrible 3-vertices. By Lemmas \ref{lem:1-thread} and \ref{lem:terrible-3-vertex}, we always have $d_{G}(w_{i})=d_{G^{*}}(w_{i})$ for $i=1,2$.
 For $i=1,2$, let $e_{i}=vv_{i}$ and $f_{i}=v_{i}w_{i}$.

Now we prove that $v$ has at least two $4^{+}$-neighbors in $G^{*}$. If not, let $v_{3}$ and $v_{4}$ be its  two $3^{-}$-neighbors in $G^{*}$.  By the minimality of $G$, the graph $G-\{v_{1},v_{2}\}$ has a good coloring $\phi$.

If $w_{1}=w_{2}$ or $d_{G}(w_{1})=d_{G}(w_{2})=2$,
 then it is easy to check that $|F_{\phi}(e_{i})|\le\Delta+7$ and 
$|\textit{SF}_{\phi}(f_{i})|\le\Delta+3$ for  $i=1,2$. 
Thus we have $|A_{\phi}(e_{i})|\ge\Delta-3=2$ and 
$|\textit{SA}_{\phi}(f_{i})|\ge\Delta+1=6$ for  $i=1,2$. 
Obviously we can extend $\phi$ to a good coloring of $G$ by greedily coloring the four edges $e_{1},e_{2},f_{1},f_{2}$, a contradiction.

Thus we now assume that $w_{1}\neq w_{2}$ and $w_{1}$ is a terrible 3-vertex. 
Now if $w_{2}$ is a 2-vertex, 
it is easy to see that $|F_{\phi}(e_{1})|\le\Delta+8$, $|F_{\phi}(e_{2})|\le\Delta+7$, $|\textit{SF}_{\phi}(f_{1})|\le\Delta+5$ and 
$|\textit{SF}_{\phi}(f_{2})|\le\Delta+3$. 
Thus we have
 $|A_{\phi}(e_{1})|\ge\Delta-4=1$, $|A_{\phi}(e_{2})|\ge\Delta-3=2$, $|\textit{SA}_{\phi}(f_{1})|\ge\Delta-1=4$ and 
$|\textit{SA}_{\phi}(f_{2})|\ge\Delta+1=6$. 
Obviously,  $\phi$ can be extended  to a good coloring  of $G$, a contradiction.

And if $w_{2}$ is also a terrible 3-vertex,  by Lemma \ref{lem:terrible-3-vertex-3-neighbor}, we must have $w_{1}w_{2}\notin E(G)$.
 For  $i=1,2$, let  $x_{i}$ be the other 2-neighbor of $w_{i}$ and let $g_{i}$ denote the edge $w_{i}x_{i}$.
 It is clear that both $x_{1}$ and $x_{2}$ are poor 2-vertices.
Because  $w_{1}$ and $w_{2}$ are  terrible 3-vertices, $x_{1}\neq x_{2}$  and $x_{1}x_{2}\notin E(G)$ as otherwise it contradicts Lemma \ref{lem:poor-2-vertex}. 
Now we erase the colors of $g_{1}$ and $g_{2}$ in $\phi$ yielding a new good partial coloring $\psi$  of $G$.
It is straightforward to check that
$|F_{\psi}(e_{i})|\le\Delta+7$,
$|\textit{SF}_{\psi}(f_{i})|\le\Delta+4 $ and  $|\textit{SF}_{\psi}(g_{i})|\le2\Delta$
for  $i=1,2$.
It follows that
$|A_{\psi}(e_{i})|\ge\Delta-3=2$,
and $|\textit{SA}_{\psi}(f_{i})|\ge\Delta=5$ and  $|\textit{SA}_{\psi}(g_{i})|\ge4$
for  $i=1,2$.
We can extend $\psi$ to a new good partial coloring $\sigma$ of $G$ by greedily coloring $e_{1},e_{2},g_{1},g_{2}$.
It is easy to see that   $|\textit{SA}_{\sigma}(f_{i})|\ge\Delta-4=1$
for  $i=1,2$. 
Notice that  $f_{1}$ and $f_{2}$ are at distance 3 in $G$, they can be colored with a common strongly available color.
We can further extend $\sigma$ to a good coloring of $G$, again a contradiction. The lemma follows.
\end{proof}

\begin{lemma}\label{lem:5-vertex-1-poor-neighbor}
Suppose that $\Delta^*\ge5$.
Let $v$ be  a  $5$-vertex  of $G^*$.  If  $v$ has exactly one poor  2-neighbor, then $v$ has at least one $4^{+}$-neighbor in $G^{*}$. 
\end{lemma}
\begin{proof}
	Let $v_{1}$  be the  poor  2-neighbor of $v$. By Observation \ref{obser:poor-2-vertex}, $\Delta=5$. Let $w_{1}$  denote the other  neighbor of $v_{1}$.
	 By the definition  of the poor 2-vertex, $w_{1}$ is either a 2-vertex or a terrible 3-vertex. 
	 Let $e_{1}=vv_{1}$ and $f_{1}=v_{1}w_{1}$.
		Now we prove that $v$ has at least one $4^{+}$-neighbors. If not, the other four neighbors of $v$ are all  $3^{-}$-vertices in $G^{*}$.
		By the minimality of $G$, the graph $G-v_{1}$ has a good coloring $\phi$.

	If $w_{1}$ is a  2-vertex, then $|F_{\phi}(e_{1})|\le13$ and 
	$|\textit{SF}_{\phi}(f_{1})|\le\Delta+4$. 
	Therefore,  $|A_{\phi}(e_{1})|\ge2\Delta-9=1$ and 
	$|\textit{SA}_{\phi}(f_{1})|\ge\Delta=5$. 
	The coloring  $\phi$ can be extended to a good coloring of $G$ by greedily coloring   $e_{1},f_{1}$, a contradiction.	
	Thus we may assume that	 $w_{1}$ is a terrible 3-vertex. Let  $x_{1}$ be the other 2-neighbor of $w_{1}$ and let $g_{1}$ denote the edge $w_{1}x_{1}$.
	Now we erase the color of $g_{1}$ in $\phi$ to obtain a new good partial coloring $\psi$  of $G$.
	It is straightforward to check that
	$|F_{\psi}(e_{1})|\le13$
	and $|\textit{SF}_{\psi}(f_{1})|\le\Delta+5 $ and  $|\textit{SF}_{\psi}(g_{1})|\le2\Delta$.
	Thus we have
	$|A_{\psi}(e_{1})|\ge2\Delta-9=1$
	and $|\textit{SA}_{\psi}(f_{1})|\ge\Delta-1=4$ and  $|\textit{SA}_{\psi}(g_{1})|\ge4$.
	We can extend  $\psi$ to a  good  coloring  of $G$ by greedily coloring $e_{1},f_{1},g_{1}$ one by one, again a contradiction. The lemma is proved.
\end{proof}

\begin{lemma}\label{lem:6+-vertex-Delta-2}
Suppose that $\Delta^*\ge6$.
Let $v$ be  a  $\Delta^*$-vertex  of $G^*$.  
If  $v$ has exactly $\Delta^*-2$ poor  2-neighbors, then $v$ has  exactly two $4^{+}$-neighbors in $G^{*}$. 
\end{lemma}
\begin{proof}
Since $v$ has poor 2-neighbors, by Observation \ref{obser:poor-2-vertex}, we have $\Delta^*=\Delta$. 	Let  $v_{1},v_{2},\ldots,v_{\Delta-2}$ be the   $\Delta-2$ poor  2-neighbors of $v$  and let $v_{\Delta-1}$ and $v_{\Delta}$ be the other two neighbors of $v$.  	
	By Observation \ref{obser:poor-neighbor-nonadjacent}, $v_{i}v_{j}\notin E(G)$ for any two distinct integers $i,j\in[1,\Delta-2]$. For each $i\in[1,\Delta-2]$, let $w_{i}$ denote the another neighbor of $v_{i}$ and let $e_{i}=vv_{i}$ and $f_{i}=v_{i}w_{i}$.
	For each $i\in[1,\Delta-2]$, it is clear that  $w_{i}$ is either a 2-vertex or a terrible 3-vertex and thus  $d_{G}(w_{i})=d_{G^{*}}(w_{i})$ according to Lemmas \ref{lem:1-thread} and \ref{lem:terrible-3-vertex}. 
	Now we prove that $d_{G^{*}}(v_{\Delta-1})\ge4$ and $d_{G^{*}}(v_{\Delta})\ge4$. If not, we may assume that   $d_{G^{*}}(v_{\Delta-1})\le3$.

	By the minimality of $G$, the graph $G-\{v_{1},v_{2},\ldots,v_{\Delta-2}\}$ has a good coloring $\phi$.  Since  $d_{G^{*}}(v_{\Delta-1})\le3$, we must have $|F_{\phi}(e_{i})|\le\Delta+5$ and 
	$|\textit{SF}_{\phi}(f_{i})|\le\Delta+4$ for each $i\in[1,\Delta-2]$. 
It follows that $|\textit{SA}_{\phi}(e_{i})|\ge\Delta-1$ and 
	$|\textit{SA}_{\phi}(f_{i})|\ge\Delta$ for each $i\in[1,\Delta-2]$. 
Therefore, we can further extend $\phi$ to a new good partial coloring $\psi$ of $G$ by greedily coloring $e_{1},e_{2},\ldots,e_{\Delta-2}$, under which 	$|\textit{SA}_{\psi}(f_{i})|\ge2$ for each $i\in[1,\Delta-2]$. 
Notice that for each $i\in[1,\Delta-2]$, there is at most one edge $f_{j}$ at distance one from $f_{i}$ in $G$, where $j\in[1,\Delta-2]$. 
Recall that	 $w_{i}$ is either a 2-vertex or a terrible 3-vertex  in $G^{*}$ for each $i\in[1,\Delta-2]$.
By Lemma \ref{lem:3-vertex}, for each $i\in[1,\Delta-2]$, there is at most one integer $j\in[1,\Delta-2]\setminus \{i\}$ such that $d_{G}(f_{i},f_{j})=1$. Let $i$ and $j$ be two distinct integers in $[1,\Delta-2]$. If $d_{G}(f_{i},f_{j})\neq1$, then due to   Lemmas \ref{lem:terrible-3-vertex-3-neighbor} and \ref{lem:poor-2-vertex},  $d_{G}(f_{i},f_{j})=3$.
Therefore, 
the coloring $\psi$ can be further extended to a good coloring of $G$ by greedily coloring  $f_{1},f_{2},\ldots,f_{\Delta-2}$. This gives a contradiction. The lemma is proved.	
\end{proof}

\begin{lemma}\label{lem:6+-vertex-Delta-3}
	Suppose that $\Delta^*\ge6$.
Let $v$ be  a  $\Delta^*$-vertex  of $G^*$. 
	If  $v$ has exactly $\Delta^*-3$ poor  2-neighbors, then $v$ has at least one $4^{+}$-neighbor in $G^{*}$. 
\end{lemma}
\begin{proof}
Since $v$ has poor 2-neighbors, by Observation \ref{obser:poor-2-vertex}, it holds that $\Delta^*=\Delta$. 		Let  $v_{1},v_{2},\ldots,v_{\Delta-3}$ be the   $\Delta-3$ poor  2-neighbors of $v$.  
	By Observation \ref{obser:poor-neighbor-nonadjacent}, $v_{i}v_{j}\notin E(G)$ for any two distinct $i,j\in[1,\Delta-3]$. 
	For each $i\in[1,\Delta-3]$, let $w_{i}$ denote the another neighbor of $v_{i}$ and let $e_{i}=vv_{i}$ and $f_{i}=v_{i}w_{i}$.
It is clear that  each $w_{i}$ is either a 2-vertex or a terrible 3-vertex and  thus $d_{G}(w_{i})=d_{G^{*}}(w_{i})$ according to Lemmas \ref{lem:1-thread} and \ref{lem:terrible-3-vertex}, where  $i\in[1,\Delta-3]$. 

Suppose on the contrary that $v$ has no  $4^{+}$-neighbors in $G^{*}$. 
By the minimality of $G$, the graph $G-\{v_{1},v_{2},\ldots,v_{\Delta-3}\}$ has a good coloring $\phi$.  It is easy to check that $|F_{\phi}(e_{i})|\le11$ and 
$|\textit{SF}_{\phi}(f_{i})|\le\Delta+5$ for each $i\in[1,\Delta-3]$. 
It follows that $|A_{\phi}(e_{i})|\ge2\Delta-7\ge\Delta-3$ and 
$|\textit{SA}_{\phi}(f_{i})|\ge\Delta-1$ for each $i\in[1,\Delta-3]$. 
Therefore, we can extend $\phi$ to a new good partial coloring $\psi$ of $G$ by greedily coloring $e_{1},e_{2},\ldots,e_{\Delta-3}$, in which 	$|\textit{SA}_{\psi}(f_{i})|\ge2$ for each $i\in[1,\Delta-3]$.

Recall that	 $w_{i}$ is either a 2-vertex or a terrible 3-vertex  in $G^{*}$ for each $i\in[1,\Delta-4]$.
By Lemma \ref{lem:3-vertex}, for each $i\in[1,\Delta-3]$, there is at most one integer $j\in[1,\Delta-3]\setminus \{i\}$ such that $d_{G}(f_{i},f_{j})=1$. Let $i$ and $j$ be two distinct integers in $[1,\Delta-3]$. If $d_{G}(f_{i},f_{j})\neq1$, then due to   Lemmas \ref{lem:terrible-3-vertex-3-neighbor} and \ref{lem:poor-2-vertex},  $d_{G}(f_{i},f_{j})=3$.
Therefore,
the coloring $\psi$ can be further extended to a good coloring of $G$ by greedily coloring  $f_{1},f_{2},\ldots,f_{\Delta-3}$, a contradiction. The lemma holds.	
\end{proof}

We have now all the ingredients to prove the case of maximum average degree less than $14/5$.
\begin{theorem}[First case of Theorem~\ref{Main-th-mad}]
Let  $G$ be a  graph with maximum degree $\Delta$. 

If $\mad(G)<14/5$, then $\chi_{ss}'(G)\le 2\Delta+4$.
\end{theorem}
\begin{proof}
Assume for contradiction that the Theorem does not hold, and 
  let   $G$ be a counterexample   with the fewest $2^{+}$-vertices, and subject to this with the fewest edges. Denote by $G^*$  the graph obtained from $G$ by deleting all its 1-vertices.
For each vertex $v\in V(G^{*})$, we use $w(v)$ and $w^{*}(v)$ to denote the  \emph{initial weight} and  the \emph{final weight}, respectively.

Let $w(v)=d_{G^{*}}(v)$ for each $v\in V(G^{*})$.
Because $G$ has maximum average degree less than ${14}/{5}$,  the average degree of  $G^{*}$ does not exceed ${14}/{5}$. 
Hence, we have 
$$\sum_{v\in V(G^{*})}w(v)=\sum_{v\in V(G^{*})}d_{G^{*}}(v)<\frac{14}{5}|V(G^{*})|.$$

The following are the discharging rules.

\noindent{\textbf{(R1)}.}  Each bad $3$-vertex gives $\frac{2}{5}$ to each of its nonpoor 2-neighbors.

\noindent{\textbf{(R2)}.}  If a good $3$-vertex $v$ has a bad 3-neighbor $u$, then we do one of the following:

\textbf{(R2-1)}. If $u$ has a bad 3-neighbor, then $v$ gives  $\frac{1}{5}$ to $u$;

\textbf{(R2-2)}. Otherwise,  $v$ gives  $\frac{1}{10}$ to $u$.

\noindent{\textbf{(R3)}.}  Each  $4^{+}$-vertex  gives $\frac{2}{5}$ to each of its  nonpoor 2-neighbors  and   $\frac{1}{5}$ to each of its  bad 3-neighbors.

\noindent{\textbf{(R4)}.}  If $G^*$ has poor 2-vertices, then  each  $\Delta^*$-vertex  gives $\frac{4}{5}$ to each of its poor 2-neighbors.

\begin{claim}\label{claim:w-14/5}
	For each $v\in V(G^{*})$, $w^{*}(v)\ge {14}/{5}$.
\end{claim}
\begin{clproof}
	According to Lemmas \ref{lem:delta} and \ref{lem:Delta}, $G^{*}$ has no 1-vertices and $\Delta^*\ge4$.

	Suppose that  $d_{G^{*}}(v)=2$.
	If $v$ is a poor 2-vertex, then according to  Lemma \ref{lem:poor-2-vertex},
	$v$ has a $\Delta$-neighbor in $G^{*}$.
	By \textbf{(R4)}, we  have
	$w^{*}(v)=w(v)+\frac{4}{5}=2+\frac{4}{5}={14}/{5}$. And if  $v$ is nonpoor, then each neighbor of $v$ is  either a  bad $3$-vertex  
	 or a $4^{+}$-vertex.  By  \textbf{(R1)} and  \textbf{(R3)}, each neighbor of $v$ gives $\frac{2}{5}$ to  $v$ and thus  $w^{*}(v)= w(v)+2\times\frac{2}{5}=2+\frac{4}{5}={14}/{5}$.
	
		Suppose that  $d_{G^{*}}(v)=3$.
	If $v$ is terrible, then $v$ does not give or receive charge and thus
	$w^{*}(v)= w(v)=3>{14}/{5}$.
	
		If $v$ is bad, $v$ has exactly one 2-neighbor, denoted as $w$, then by Lemma \ref{lem:bad-3-vertex-bad-neighbor}, $w$ is nonpoor and  $v$ has  at most one bad 3-neighbor and has no terrible 3-neighbors.
		 By  \textbf{(R1)}, $v$ needs to give $\frac{2}{5}$ to its nonpoor 2-neighbor $w$.  Now if $v$
		has exactly one bad 3-neighbor, then the third neighbor of  $v$  is either a good 3-vertex (see Lemma	\ref{lem:bad-3-vertex}) or a 
	$4^{+}$-vertex, by \textbf{(R2-1)} or \textbf{(R3)}, this neighbor will give   $\frac{1}{5}$ to $v$, and thus  $w^{*}(v)= w(v)-\frac{2}{5}+\frac{1}{5}=3-\frac{1}{5}={14}/{5}$.  
		And if $v$
		has no bad 3-neighbors, then each $3^{+}$-neighbor of $v$ is either a  good 3-vertex or  a $4^{+}$-vertex,  by \textbf{(R2-2)} or \textbf{(R3)}, each  $3^{+}$-neighbor  gives at least  $\frac{1}{10}$ to $v$, and thus  $w^{*}(v)\ge w(v)-\frac{2}{5}+2\times\frac{1}{10}=3-\frac{1}{5}={14}/{5}$.

	If	$v$ is a good 3-vertex, then $v$ has no 2-neighbors. 
	If $v$ has no bad 3-neighbors, then  $v$ does not give or receive charge and thus
	$w^{*}(v)= w(v)=3>{14}/{5}$.
	We may assume that   $v$ has  bad 3-neighbors.
  By Lemma \ref{lem:3-vertex-bad-neighbor}, $v$ has at most two bad 3-neighbors.	
  If $v$ has exactly two bad 3-neighbors, then by Lemma \ref{lem:bad-3-vertex},  the two bad 3-neighbors of $v$ have no bad 3-neighbors.
  Then  by \textbf{(R2-2)}, $v$ gives $\frac{1}{10}$ to each of its   bad 3-neighbors  and thus    $w^{*}(v)= w(v)-2\times\frac{1}{10}=3-\frac{1}{5}={14}/{5}$.
Now we assume that  $v$ has exactly one bad 3-neighbor, denoted by $u$. 
	If $u$  has a bad 3-neighbor, then by Lemma \ref{lem:bad-3-vertex}, the other two neighbors of $v$ are all $4^{+}$-vertices in $G^{*}$,
	 by \textbf{(R2-1)}, $v$ gives $\frac{1}{5}$ to $u$ and thus    $w^{*}(v)= w(v)-\frac{1}{5}=3-\frac{1}{5}={14}/{5}$.
	And if $u$  has no bad 3-neighbor, then 
	 by \textbf{(R2-2)}, $v$ gives $\frac{1}{10}$ to  $u$  and thus    $w^{*}(v)= w(v)-\frac{1}{10}=3-\frac{1}{10}=\frac{29}{10}>{14}/{5}$.

The remaining proof are divided into three cases according to the value of $\Delta^*$.

\textbf{Case 1.} $\Delta^*=4$.

Let $v$ be a 4-vertex in $G^{*}$. By Lemma \ref{lem:4-vertex-poor-neighbor}, $v$ has at most one poor 2-neighbor. If $v$   has exactly one poor 2-neighbor, again by Lemma \ref{lem:4-vertex-poor-neighbor}, $v$ has at least two 4-neighbors and thus $v$ has at most one neighbor that is a nonpoor 2-vertex or a bad 3-neighbor. By \textbf{(R3)} and \textbf{(R4)}, we must have 
  $w^{*}(v)\ge w(v)-\frac{2}{5}-\frac{4}{5}=4-\frac{6}{5}={14}/{5}$.

And if $v$   has no poor 2-neighbors, by Lemma \ref{lem:4-vertex}, $v$ has at most three neighbors that are nonpoor 2-vertices or bad 3-vertices. Then by \textbf{(R3)},  $w^{*}(v)\ge w(v)-3\times\frac{2}{5}=4-\frac{6}{5}={14}/{5}$.

\textbf{Case 2.} $\Delta^*=5$.

Let $v$ be a 4-vertex in $G^{*}$.
By Lemma \ref{lem:poor-2-vertex}, $v$ has no poor 2-vertices.
By Lemma  \ref{lem:4-vertex},  $v$ has 
 at most three neighbors that are nonpoor 2-vertices or bad 3-vertices.
 Thus by \textbf{(R3)},  $w^{*}(v)\ge w(v)-3\times\frac{2}{5}=4-\frac{6}{5}={14}/{5}$. 

Let $v$ be a 5-vertex in $G^{*}$.  By Lemma \ref{lem:5-vertex}, $v$ has at most two poor 2-neighbors. 
If $v$ has exactly  two poor 2-neighbors, by Lemma \ref{lem:5-vertex-2-poor-neighbor}, $v$ has at least two $4^{+}$-neighbors in $G^{*}$ and thus at most one neighbor of $v$ is a nonpoor 2-vertex or a bad 3-vertex, by \textbf{(R3)} and \textbf{(R4)}, we must have 
$w^{*}(v)\ge w(v)-\frac{2}{5}-2\times\frac{4}{5}=5-2=3>{14}/{5}$. 

If $v$ has exactly one  poor 2-neighbor, by Lemma \ref{lem:5-vertex-1-poor-neighbor}, $v$ has at least one $4^{+}$-neighbors in $G^{*}$ and thus $v$ has  at most three neighbors that are nonpoor 2-vertices or bad 3-vertices, by \textbf{(R3)} and \textbf{(R4)}, we must have 
$w^{*}(v)\ge w(v)-3\times\frac{2}{5}-\frac{4}{5}=5-2=3>{14}/{5}$. 

If $v$ has no  poor 2-neighbors, then $v$  at most five neighbors that are nonpoor 2-vertices or  bad 3-vertices, by \textbf{(R3)} , we must have 
$w^{*}(v)\ge w(v)-5\times\frac{2}{5}=5-2=3>{14}/{5}$. 

\textbf{Case 3.} $\Delta^*\ge6$.

Let $v$ be a 4-vertex in $G^{*}$. By Lemma \ref{lem:poor-2-vertex}, $v$ has no poor 2-vertices.
By Lemma  \ref{lem:4-vertex},  $v$ has 
at most three neighbors that are nonpoor 2-vertices or bad 3-vertices.  
 Thus by \textbf{(R3)},  $w^{*}(v)\ge w(v)-3\times\frac{2}{5}=4-\frac{6}{5}={14}/{5}$. 

Suppose that $5\le d_{G^{*}}(v) \le \Delta^*-1$. By Lemma \ref{lem:poor-2-vertex}, $v$ has no poor 2-vertices.
 By 
\textbf{(R3)}, we must have $w^{*}(v)\ge w(v)-d_{G^{*}}(v)\times\frac{2}{5}=d_{G^{*}}(v)-d_{G^{*}}(v)\times\frac{2}{5}=\frac{3}{5}\times d_{G^{*}}(v)\ge\frac{3}{5}\times 5=3 >{14}/{5}$. 

Suppose that $ d_{G^{*}}(v) = \Delta^*$. By Lemma \ref{lem:Delta-vertex-poor-neighbor}, $v$ has at most $\Delta^*-2$ poor 2-neighbors. Now if $v$ has exactly  $\Delta^*-2$ poor 2-neighbors, by Lemma \ref{lem:6+-vertex-Delta-2}, the other two neighbors of $v$ are all $4^{+}$-neighbors in $G^{*}$. Then by \textbf{(R4)}, we must have  $w^{*}(v)= w(v)-(\Delta^*-2)\times\frac{4}{5}=\Delta^*-(\Delta^*-2)\times\frac{4}{5}=\frac{1}{5}\times \Delta^*+\frac{8}{5}\ge\frac{1}{5}\times6+\frac{8}{5} ={14}/{5}$. 

And if $v$ has exactly  $\Delta^*-3$ poor 2-neighbors, by Lemma \ref{lem:6+-vertex-Delta-3}, $v$ has at least one  $4^{+}$-neighbor in $G^{*}$ and thus  $v$ has at most two neighbors that are nonpoor 2-vertices or bad 3-neighbors.
 Then by \textbf{(R3)} and \textbf{(R4)}, we must have  $w^{*}(v)\ge w(v)-2\times\frac{2}{5}-(\Delta^*-3)\times\frac{4}{5}=\frac{1}{5}\times \Delta^*+\frac{8}{5}\ge\frac{1}{5}\times6+\frac{8}{5} ={14}/{5}$. 
 
 Finally,  if $v$ has at most  $\Delta^*-4$ poor 2-neighbors,  we may assume that $v$ has exactly $i$ ($0\le i \le\Delta^*-4$) poor 2-neighbors and thus   $v$ has at most $\Delta^*-i$ neighbors that are nonpoor 2-vertices or bad 3-neighbors.
 Then by \textbf{(R3)} and \textbf{(R4)}, we must have  $w^{*}(v)\ge w(v)-(\Delta^*-i)\times\frac{2}{5}-i\times\frac{4}{5}=\frac{3}{5}\times \Delta^* -\frac{2}{5}\times i
 \ge\frac{3}{5}\times \Delta^*-\frac{2}{5}\times(\Delta^*-4)
 =\frac{1}{5}\times \Delta^* +\frac{8}{5} \ge
 \frac{1}{5}\times6+\frac{8}{5} ={14}/{5}$. 

We have proved that $w^{*}(v)\ge{14}/{5}$
for any $v\in V(G^{*})$. 
	The claim is true.
\end{clproof}

It follows from Claim \ref{claim:w-14/5} that
$$\frac{14}{5}|V(G^{*})|\le \sum_{v\in V(G^{*})}w^{*}(v)=\sum_{v\in V(G^{*})}w(v)<\frac{14}{5}|V(G^{*})|,$$
which is a contradiction. 
\end{proof}

\section{Discussion}

Recall that  Lu{\v{z}}ar, Mockov{\v{c}}iakov{\'a} and Sot{\'a}k \cite{LMS2024} proposed the following conjecture.
\LMSConj*
If this conjecture is true, then the upper bound $2\Delta+4$ is the best possible due to Example~\ref{example}.
However,  all the examples in  Example~\ref{example} are planar graphs with  girth  3. It is not clear whether there exist planar graphs with  girth greater than 3 and less than 8 and whose semistrong chromatic indices are exactly $2\Delta+4$. 
We further propose the following.

\begin{problem}\label{Pro:no3cycle}
Is it true that $\chi_{ss}'(G)\le 2\Delta+2$ for	every triangle-free planar graph $G$ with maximum degree $
\Delta$?
\end{problem}

If the answer to \Cref{Pro:no3cycle} is affirmative, then the 5-prism (i.e., a 3-regular planar graph with girth 4) shows the sharpness of the upper bound $2\Delta+2$.
Intuitively, the semistrong chromatic index of a planar graph relies on its maximum degree and girth.

\begin{problem}
    Give a precise upper bound for $\chi_{ss}'(G)$ of a planar graph $G$ in terms of $\Delta$ and $\g(G)$.
\end{problem}

Recall that for every graph $G$ we have $\chi'(G)\le a'(G) \le \chi_{ur}'(G)\le \chi_{ss}'(G) \le \chi_{s}'(G)$.
This motivates us to propose an analogous version of \Cref{conj-planar-semi-L}, where the semistrong chromatic index is replaced by the uniquely restricted chromatic index.

\begin{conjecture}  \label{our-conj-planar-ur}
Every  planar graph $G$  with maximum degree $\Delta$ satisfies that $\chi_{ur}'(G)\le 2\Delta+4$.
\end{conjecture}

For a planar graph $G$, Example \ref{example} with $\Delta=4$ shows that $2\Delta+4$ would be the best possible upper bound for  $\chi_{ur}'(G)$.\\

An application of strong edge coloring is to model the conflict-free channel assignment problem in radio networks \cite{R1997,NKGB2000}.
However, when the channel resource is limited, it is possible that  one cannot construct a proper strong edge coloring.
In such a case, a certain kind of relaxation becomes necessary when assigning channels to transmitters.
This may be one of the reasons why the various relaxed strong edge colorings have attracted attention in recent years, refer to  \cite{LMS2024} for the classification of various relaxations of strong edge coloring.
Admittedly, it remains unclear whether relaxations such as the semistrong-form  and uniquely-restricted-form  correspond directly to specific patterns in real-world applications. Nevertheless, exploring their potential applications and understanding the comparative advantages of different edge coloring models is a meaningful direction for future work.

Interestingly, in the  planar case, these two relaxations may reduce the chromatic number by nearly half (refer to Table \ref{tab:bound}).
It encourages us to explore whether these two  chromatic indices can also be much lower than the strong  chromatic index on other special classes of graphs.

Note that  every graph  $G$ satisfies $\chi_{ur}'(G)\le\Delta^{2}$, with equality  if and only if $G$ is the complete bipartite graph $K_{\Delta,\Delta}$ \cite{BRS2019}.
Surprisingly, despite the stronger constraint of the semistrong edge coloring, it yields the same result, meaning that $\chi_{ss}'(G)\le\Delta^{2}$ holds for every graph $G$ and the equality holds  if and only if $G$ is isomorphic to $K_{\Delta,\Delta}$ \cite{LMS2024,LL2023}.
If \Cref{conj-planar-semi-L} is proven,  it would follow that the two upper bounds of $\chi_{ur}'(G)$ and $\chi_{ss}'(G)$
are also asymptotically equal (up to a constant) for planar graphs.
Then, a natural problem is whether these two chromatic indices also maintain a similar relationship on other special classes of graphs.\\

\noindent{\bf Acknowledgements}
We are sincerely grateful to Patrice Ossona de Mendez for his invaluable suggestions, constructive revisions and helpful discussions throughout the writing of this article. 
We would also like to thank the anonymous reviewers for their careful reading and valuable comments, which have helped improve the quality of this work.\\

\noindent{\bf Funding}  The first author is supported by  SEU Innovation Capability Enhancement Plan for Doctoral Students (CXJH\_SEU 24119).
The second author is supported by NSFC 11771080.\\

\noindent{\bf  Data availability}
Data sharing not applicable to this article as no datasets were generated or analysed during
the current study.\\

\section*{Declarations}

\noindent{\bf Conflict of interest} The authors declare that they have no conflict of interest.

\bibliographystyle{amsplain}
\bibliography{ref}

\end{document}